\documentclass[12pt]{amsart}

\usepackage{amsmath,amsthm,amssymb,color,array}

\newtheorem{theorem}{Theorem}[section]
\newtheorem{conjecture}[theorem]{Conjecture}
\newtheorem{prop}[theorem]{Proposition}
\newtheorem{lemma}[theorem]{Lemma}

\newtheorem{cor}[theorem]{Corollary}
\theoremstyle{remark}
\newtheorem{remark}[theorem]{Remark}

\theoremstyle{definition}
\newtheorem{definition}[theorem]{Definition}

\numberwithin{equation}{section}
\newcommand{\lbl}[1]{\label{#1}}

\begin{document}

\title{Abelian Functions for Trigonal Curves of Genus Three}
\author{J. C. Eilbeck}
\address{Department of Mathematics and the Maxwell Institute for Mathematical 
Sciences, Heriot-Watt University, Edinburgh,
   UK EH14 4AS}
\email{J.C.Eilbeck@hw.ac.uk}
\author{V. Z. Enolski} 
\address{Institute of Magnetism, Vernadski blvd. 36, Kiev-142, 
Ukraine}
\email{vze@ma.hw.ac.uk}

\author{S. Matsutani} 
\address{8-21-1 Higashi-Linkan, Sagamihara, 228-0811, Japan}
\email{rxb01142@nifty.com}

\author{Y. \^Onishi}
\address{Faculty of Humanities and Social Sciences
Iwate University, Ueda 3-18-34, Morioka 020-8550, Japan}
\email{onishi@iwate-u.ac.jp}

\author{E. Previato}
\address{Department of Mathematics and Statistics,
Boston University,
Boston, MA 02215-2411,
USA}
\email{ep@bu.edu}

\begin{abstract}
We develop the theory of generalized Weierstrass $\sigma$- and $\wp$-functions 
defined on a trigonal curve of genus three.  In particular
we give a list of the associated partial differential equations
satisfied by the $\wp$-functions, a proof that the coefficients of 
the power series expansion of the $\sigma$-function are polynomials of moduli 
parameters, and the derivation of two addition formulae.
\end{abstract}

\maketitle

\section*{Introduction} 

Constructive theories of Abelian and modular functions associated with
algebraic curves have seen an upsurge of interest in recent times.
These classical functions have been of crucial importance in
mathematics since their definition at the hands of Abel, Jacobi,
Poincar\'e and Riemann, but their relevance in physics and applied
mathematics has greatly developed over the past three decades.
Algebraic curves are here intended as Riemann surfaces, unless
specified to be singular.

The study of the simplest hyperelliptic curves, namely curves of genus
two, goes back to the beginning of the 20th century, and these are
treated in much detail in advanced textbooks, see for example Baker
(1907) \cite{ba07} and more recently Cassels and Flynn (1996)
\cite{cafl96}.  Not so much is known about the simplest trigonal
curves, which have genus three.  The study of modular functions of
these curves was originated by Picard, and reprised recently by Shiga
\cite{shiga88} and his school.  In this paper we study Abelian
functions associated with the simplest general type of curve, the
general (3,4) curve.   
% by analyzing a special
% case, namely a ``purely trigonal'' curve\footnote{Some papers call a
%   curve of this type a ``$\mathbb{Z}_3$-curve'', where $\mathbb{Z}_3$
%   denotes $\mathbb{Z}/3\mathbb{Z}$.}  in the terminology introduced in
% \cite{on05}. 
This is an $(n,m)$-curve in the sense of Burchnall-Chaundy \cite{bc28}
%  with $(n,m)=(3,4)$; it is also a cyclic
% Galois cover of $\mathbb{P}^1$, of genus three, so it is not a general
% curve of genus three.

Our work is based on the realization of Abelian functions as
logarithmic derivatives of the multi-dimensional $\sigma$-function.
This approach is due to Weierstrass and Klein and was developed by
Baker \cite{ba97}; for recent developments of the theory of
multi-dimensional $\sigma$-functions, see Grant \cite{gr90},
Buchstaber, Enolskii, and Leykin \cite{bel97}, 
Buchstaber and Leykin \cite{bl02}, \cite{bl05}, 
Eilbeck, Enolskii and Previato \cite{eep03}, 
and Baldwin and Gibbons \cite{bg06}, among others.  

We shall adopt as a template the Weierstrass theory of elliptic
functions, trying to extend as far as possible these results to the
case of the trigonal genus-three  curve.  Let $\sigma(u)$ and $\wp(u)$
be the standard functions in Weierstrass elliptic function theory.
They satisfy the well-known formulae
\begin{equation}
  \wp(u) = - \frac{d^2}{du^2}\log \sigma (u), 
  \quad (\wp')^2=4\wp^3-g_2\wp-g_3,
  \quad \wp''=6\wp^2-\tfrac12 g_2 \lbl{WP}
\end{equation}
and the addition formula, which is a  basic formula of the theory
\begin{equation}
   -\frac{\sigma(u+v)\sigma(u-v)}{\sigma(u)^2\sigma(v)^2}=\wp(u)-\wp(v).
   \lbl{eq0.1}
\end{equation}
 
We present here two addition formulae (Theorems \ref{T7.1} and
\ref{T8.1}).  The first of these is for the general trigonal curve of
degree four, whereas the second is restricted to a ``purely trigonal''
curve of degree four (see (\ref{eq5.1})).  The first main Theorem
\ref{T7.1} is the natural generalization of (\ref{eq0.1}).  The
authors realized the existence of the formula of the second main
Theorem \ref{T8.1} from \cite{on05}.  However we were not able to use
that paper to establish our result, instead working from results by
Cho and Nakayashiki \cite{cn06}, Grant's paper \cite{gr90}, p.100,
(1.6), or a calculation using \cite{bel00}.  The crucial part is to
identify the coefficients of the right hand sides of these two
formulae.  To calculate these, we used a power-series expansion of the
$\sigma$-function, stimulated by the works of Buchstaber and Leykin
\cite{bl02} for hyperelliptic case and of Baldwin and Gibbons
\cite{bg06} for a purely trigonal curve of genus four.

The $\sigma$-functional realization of Abelian functions of trigonal
curve of arbitrary genus $g$ was previously developed in \cite{bel00}
and \cite{eel00}. Using these results in the case of $g=3$ we present
explicit formulae for 6 canonical meromorphic differentials and the
symmetric bi-differential which allow us to derive a complete set of
relations for trigonal $\wp$-functions, generalizing the above
relations for the Weierstrass $\wp$-function.

We note that we have recently developed a parallel, but more
limited theory, for purely trigonal curves of genus {\em four} in
\cite{bego06}, a paper which draws heavily on the results presented
here.  It is perhaps useful to compare and contrast these two cases.
As demonstrated in Schilling's generalization of the Neumann system
\cite{sc89}, there are basically two cases of trigonal cyclic covers,
the order of a related linear differential operator that commutes with
the given one of order three being congruent to 1 or 2 modulo 3,
respectively.  In each case, the action variables of the integrable
system parametrize a family of curves of the same type, thus the
family of curves in the $(3,4)$-case cannot be obtained as a limit of
that in the $(3,5)$-case, as they have different dimensions.  In the
present paper, we develop the method and prove the addition formulae,
together with the characterising differential equations, for the
former case, in that the highest power of $x$ appearing in the
equation of the curve is 4 ($\equiv 1$ modulo 3); this corresponds to
the `base' case of the Boussinesq equation, the smallest-genus
spectral curve of an algebro-geometric third-order operator.  In
\cite{bego06}, the case where the highest power of $x$ appearing in
the equation of the curve is 5 ($\equiv 2$ modulo 3) is addressed.
The differences in the two cases manifest themselves in a number of
ways, for example the parity of the $\sigma$ function is different in
the two cases, and the two-term addition formulae are antisymmetric in
the genus 3 case and symmetric in the genus 4 case.  Also the results
are given for the {\em general} $(3,4)$-curve here, whereas only for the
{\em purely} trigonal $(3,5)$-case in \cite{bego06}.  It may be
possible with some work to relate the $(3,5)$-case to the $(3,4)$-case,
but this would not be straightforward and we have not yet attempted
this.
   
Our study is far from complete, and a number of questions still
remain. One of the first problems still to be considered should be the
explicit recursive construction of the $\sigma$-series generalizing
the one given by Weierstrass; for a hyperelliptic curve of genus two,
this result was found by Buchstaber and Leykin \cite{bl02}, who also
devised a procedure to derive such recursions for the whole family of
$(n,m)$-curves \cite{bl02}, \cite{bl05}. Another problem is the deeper
understanding of the algebraic structure of the addition theorems
developed here, in order to generalize results to higher genera.  As a
pattern one can consider the addition formula of \cite{bel97} for
hyperelliptic $\sigma$-functions of arbitrary genera written in terms
of certain Pfaffians.  Also, the description of Jacobi and Kummer
varieties as projective varieties, whose coordinates are given in
terms of (derivatives of) trigonal $\wp$-functions, is far from
complete.  We hope the results we present to be the first steps
towards a general theory of trigonal curves of arbitrary genus, as
well as a tool in the study of projective varieties which are images
of Jacobians.

The paper is organized as follows.  We first discuss the basic
properties of the general $(3,4)$-curve in Section \ref{BasicProps},
and define a restricted version of this curve, the ``purely trigonal
case'', in Section \ref{PureTrig}.  In Section \ref{sigma}, we
introduce the $\sigma$ function for the general curve, and in Section
\ref{Abelian} the Abelian functions $\wp_{ij}$ and their derivatives.
Section \ref{PDEs} of the paper is devoted to the various differential
relations satisfied by these Abelian functions, and the series
expansion of the $\sigma$ function is discussed in Section
\ref{sigma_expan}, in which the result (Theorem \ref{L6.1}) is new, is
proved quite constructively, and is the key for the rest of papers.
Let $\Theta^{[2]}$ be the standard theta divisor, namely the image of
the Abelian map of the symmetric square of the curve that we consider,
in its Jacobian variety $J$.  The basis of the spaces
$\varGamma(J,\mathcal{O}(n\Theta^{[2]}))$ of functions on $J$ whose
poles are at most of order $n$ along $\Theta^{[2]}$ are discussed in
Section \ref{BasisG}, as a preliminary to the two main addition
Theorems in Sections \ref{Add1} and \ref{Add2}, respectively.  The
first addition theorem is a two-term relation for the general
$(3,4)$-curve, and the second a three-term relation for the ``purely
trigonal'' $(3,4)$-curve.  Appendix A has some formulae for the
fundamental bi-differential, and Appendix B has a list of quadratic
three-index relations for the ``purely trigonal'' case only, as the
full relations would require too much space.  The web site
\cite{Weier34} contains more details of the relations omitted through
lack of space.

While Sections 1 and 2 overlap somewhat with material in \cite{on05}, 
we believe that the results are useful to make the present paper reasonably
self-contained.

\newpage
\tableofcontents
%\newpage

\section{Trigonal curves of genus three} %%%%%%%% Section 2 %%%%%%%%%%%%%%%%%
\setcounter{equation}{0}
\label{BasicProps}
Let  $C$  be the curve defined by  $f(x,y)=0$, where 
\begin{align}
\begin{split}
   f(x,y) =
   y^3 + &(\mu_1 x + \mu_4)y^2 + (\mu_2 x^2 + \mu_5 x + \mu_8)y \\
    & - (x^4 + \mu_3x^3 + \mu_6x^2 + \mu_9x + \mu_{12}),
    \quad \text{($\mu_j$ are constants),}%\nonumber
\end{split}\lbl{eq1.1}
\end{align}
with the unique point $\infty$ at infinity.  This curve is of genus
$3$, if it is non-singular.  We consider the set of 1-forms
\begin{equation}
   \omega_1=\frac{ dx}{f_y(x,y)},  \quad
   \omega_2=\frac{xdx}{f_y(x,y)},  \quad
   \omega_3=\frac{ydx}{f_y(x,y)},
   \lbl{eq1.2} %(2.2)
\end{equation}
where  $f_y(x,y)=\frac{\partial}{\partial y}f(x,y)$.  
This is a basis of the space of differentials of the first kind on $C$.
We denote the vector consisting of the forms (\ref{eq1.2}) by
\begin{equation}
  \omega=(\omega_1,\omega_2,\omega_3)
  \lbl{eq1.2.5}
\end{equation}
We know, by the general theory, that for three variable points
$(x_1,y_1)$, $(x_2,y_2)$, and $(x_3,y_3)$ on $C$, the sum of integrals
from $\infty$ to these three points
\begin{equation}
\begin{aligned}
   u &=(u_1,u_2,u_3) \\
   & =\int_{\infty}^{(x_1,y_1)}\omega+\int_{\infty}^{(x_2,y_2)}\omega
     +\int_{\infty}^{(x_3,y_3)}\omega
\end{aligned} \lbl{eq1.3}%   \nonumber
 % \\
\end{equation}
fills the whole space $\mathbb{C}^3$.  We denote the points in
$\mathbb{C}^3$ by $u$ and $v$ etc., and their natural coordinates in
$\mathbb{C}^3$ by the subscripts $(u_1,u_2,u_3)$, $(v_1,v_2,v_3)$.  We
denote the lattice generated by the integrals of the basis
(\ref{eq1.2}) along any closed paths on $C$ by $\Lambda$.  We denote
the manifold $\mathbb{C}^3/\Lambda$, by $J$, the Jacobian variety over
$\mathbb{C}$ of $C$. We denote by $\kappa$
the natural map to the quotient group,
\begin{equation}
   \kappa:\mathbb{C}^3\rightarrow \mathbb{C}^3/\Lambda=J.
   \lbl{eq1.4}
\end{equation}
%We have  $\Lambda=\kappa^{-1}\big((0,0,0)\big)$. 
$\Lambda$ is a lattice of the space $\mathbb{C}^3$ generated by the
integrals $\oint\omega$ along any loop on $C$.  We define for $k=1$,
$2$, $3$, $\dots$, the map
\begin{equation}
\begin{aligned}
  \iota : \text{Sym}^k(C)&\rightarrow J \\
  (P_1,\cdots,P_k) &\mapsto \left(\int_{\infty}^{P_1}\omega+\cdots+
  \int_{\infty}^{P_k}\omega\right) \text{ mod}\,\Lambda,
   \end{aligned}\lbl{eq1.5}
\end{equation}
and denote its image by $W^{[k]}$.  ($W^{[k]} = J$ for $k\ge3$ by 
the Abel-Jacobi theorem.)  We will use the same symbol $u = (u_1, u_2,
u_3)$ for a point $u \in \mathbb{C}^3$ in 
$\kappa^{-1}(W^{[k]})$. 
 Let
\begin{equation}
   [-1](u_1,u_2,u_3)=(-u_1,-u_2,-u_3),
   \lbl{eq1.6}
\end{equation}
and 
\begin{equation}
   \Theta^{[k]}:=W^{[k]}\cup[-1]W^{[k]}.
   \lbl{eq1.7}
\end{equation} 
We call this  $\Theta^{[k]}$  the  $k$-th {\it standard theta subset}. 
In particular, if $k=1$, then (\ref{eq1.5}) gives an embedding of  $C$:
\begin{equation}
   \begin{aligned}
   \iota : &C\rightarrow J \\
   & P \mapsto \int_{\infty}^P\omega  \text{ mod } \Lambda.
  \end{aligned}
   \lbl{eq1.8}
\end{equation}
We note that  
\begin{equation}
   \Theta^{[2]}= W^{[2]}, \quad\Theta^{[1]}\neq W^{[1]},
\end{equation}  
differing from  the genus-3 hyperelliptic case in a suitable normalization
\cite{bel97}.  
If  $u=(u_1,u_2,u_3)$ varies on the inverse image 
$\kappa^{-1}\iota(C)=\kappa^{-1}(W^{[1]})$ of the embedded curve, 
we can take $u_3$  as a local parameter at the origin $(0,0,0)$. 
Then we have (see \cite{on05}, e.g.) 
Laurent expansions with respect to $u_3$ as follows:
\begin{equation}
   u_1=\tfrac15{u_3}^5+\cdots,\quad u_2=\tfrac12{u_3}^2+\cdots
  \lbl{u3expansion1}
\end{equation}
and 
\begin{equation}
   x(u)=\frac1{{u_3}^3}+\cdots, \quad y(u)=\frac1{{u_3}^4}+\cdots.
  \lbl{u3expansion2}
\end{equation}

We introduce a weight for several variables as follows:
\begin{definition}
  We define a {\rm weight} for constants and variables appearing in
  our relations as follows.  The weights of the variables $u_1$,
  $u_2$, $u_3$ for every $u=(u_1, u_2, u_3)$ of $W^{[k]},
  (k=1,2,\dots)$ are $5$, $2$, $1$, respectively, and the weight of
  each coefficient $\mu_j$ in {\rm(\ref{eq1.1})} is $-j$, the weights of
  $x$ and $y$ of each point $(x, y)$ of $C$ are $-3$ and $-4$,
  respectively. So, the weights of the variables are nothing but the
  order of zero at $\infty$, while the weight assigned to the
  coefficients is a device to render $f(x,y)$ homogeneous.  This is
  the reason why $\mu_7$, $\mu_{10}$, $\mu_{11}$ are absent.
\end{definition}
We remark that the weights of the variables $u_k$  are precisely 
the Weierstrass gap numbers of the Weierstrass gap sequence at $\infty$, 
whilst the weights of monomials of  $x(u)$  and  $y(u)$  correspond to 
the Weierstrass non-gap numbers in the sequence.  
In particular, in the case considered the Weierstrass gap
sequence is of the form
\[  \overline{0},\;1,\; 2,\;\overline{ 3,\; 4},\; 5,\; 
\overline{6,\; 7,\; 8,\; 9,\;10,\ldots}
\]
where orders of existing functions of the form $x^{p}y^{q}$, $p,q\in
\mathbb{N}\cup \{0\}$ are overlined.

The definition above is compatible, for instance, with the Laurent
expansion of $x(u)$ and $y(u)$ with respect to $u_3$, etc.\ for $u \in
W^{[1]}$.  Moreover, all the equalities in this paper are homogeneous
with respect to this weight.

In the next section, we use the discriminant of $C$.  Axiomatically,
the discriminant $D$ of $C$ is defined as (one of) the simplest
polynomial(s) in the $\mu_j$'s such that $D=0$ if and only if $C$ has
a singular point.  Here we are regarding $C$ as a family of curves
over $\mathbb{Z}$.  While no concrete expression of the discriminant
is necessary for the main results in this paper, we put forward a
conjecture based on the results of experimentation on special cases of
$C$ using computer algebra.
\begin{conjecture}
  Let \def\rslt{\mathrm{rslt}}
\begin{equation}
\begin{aligned}
R_1&=\rslt_x\big(\rslt_y\big(f(x,y), f_x(x,y)\big),
                 \rslt_y\big(f(x,y), f_y(x,y)\big)\big), \\
R_2&=\rslt_y\big(\rslt_x\big(f(x,y), f_x(x,y)\big), 
                 \rslt_x\big(f(x,y), f_y(x,y)\big)\big), \\
R_3&=\gcd(R_1,R_2),
\end{aligned}
\end{equation}
where $\rslt_z$ represents the resultant, namely, the determinant of
the Sylvester matrix with respect to the variable $z$.
Then $R_3$ is of weight $144$ and a perfect square 
in the ring
\begin{equation*}
\mathbb{Z}[\mu_1,\mu_4,\mu_2,\mu_5,\mu_8,\mu_3,\mu_6,\mu_9,\mu_{12}]. 
\end{equation*}
\end{conjecture}  

Unfortunately checking this condition directly is a computing task
presenting considerable difficulties due to the size of the
intermediate expressions involved.  We leave this as a conjecture and
remark only that work on a full calculation is continuing.  This
result is not crucial to this paper, but we will adopt it as a working
hypothesis (see Remark \ref{varepsilon}).  To continue, we define here
the {\it discriminant} $D$ of $C$ by a square root of $R_3$:
\begin{equation}
D=\sqrt{R_3}.
\lbl{discriminant}
\end{equation}
We comment on the choice of this root in \ref{varepsilon}.   
If the conjecture is true, $D$ is of weight  $72$. 
For the convenience of the reader we give  $R_3$
  \[
  \begin{aligned}
   R_3=\big(&256{\mu_{12}}^3-27{\mu_{12}}^2{\mu_3}^4-128{\mu_{12}}^2{\mu_6}^2
  +144{\mu_{12}}^2{\mu_6}{\mu_3}^2-192{\mu_{12}}^2{\mu_9}{\mu_3}
  +16{\mu_{12}}{\mu_6}^4\\
  &-80{\mu_{12}}{\mu_9}{\mu_6}^2{\mu_3}
  -4{\mu_{12}}{\mu_3}^2{\mu_6}^3+18{\mu_{12}}{\mu_9}{\mu_3}^3{\mu_6}
  +144{\mu_{12}}{\mu_9}^2{\mu_6}-6{\mu_{12}}{\mu_9}^2{\mu_3}^2 \\
  &-4{\mu_9}^2{\mu_6}^3-4{\mu_9}^3{\mu_3}^3
  +{\mu_9}^2{\mu_3}^2{\mu_6}^2
  +18{\mu_9}^3{\mu_6}{\mu_3}-27{\mu_9}^4\big)^6
  \end{aligned}
  \]
for the special case of $\mu_1=\mu_2=\mu_4=\mu_5=\mu_8=0$ (see Section
\ref{PureTrig}).
 
\begin{definition}
The $2$-form $\Omega((x,y),(z,w))$ on $C\times C$ is called 
{\rm fundamental 2-from of the second kind} or 
{\rm (fundamental second kind bi-differential)} if it is symmetric, namely, 
\begin{equation}
\Omega((x,y),(z,w))=\Omega((z,w),(x,y)), 
\lbl{eq3.1.6}  
\end{equation}
it has its only pole {\rm(}of second order{\rm)} along the diagonal of  $C\times C$, 
and in the vicinity of each point $(x,y)$ it is expanded in power series as 
\begin{equation}
\Omega((x,y),(z,w))=\left(\frac{1}{(\xi-\xi')^2  } +O(1)\right)d\xi d\xi'
\quad (\text{as }  (x,y)\rightarrow (z,w)), 
\lbl{expansion}
\end{equation}
where $\xi$ and $\xi'$ are local coordinates of points $(x,y)$ and $(z,w)$.  
\end{definition}
We shall look for a realization of $\Omega((x,y),(z,w))$ in the form 
\begin{equation} 
  \Omega((x,y),(z,w))=\frac{F((x,y),(z,w))dx dz}
{(x-z)^2 f_y(x,y)f_w(z,w)},\lbl{realization}   
\end{equation}
where $F((x,y),(z,w))$ is a polynomial in its variables.

\begin{lemma} {\rm (Fundamental 2-form of the second kind)} \quad % 3.1 
Let $\Sigma\big((x,y),(z,w)\big)$ be the meromorphic function on $C\times C$,
\begin{equation}
   \Sigma\big((x,y),(z,w)\big)
   =\frac{1}{(x-z)f_y(x,y)}
    \sum_{k=1}^3y^{3-k}\bigg[\frac{f(Z,W)}{W^{3-k+1}}\bigg]_W
\bigg|_{(Z,W)=(z,w)}
   \lbl{eq2.3}, % (3.3)
\end{equation} 
where $[\quad ]_W$ means removing the terms of negative powers with
respect to $W$.  Then there exist differentials $\eta_j=\eta_j(x,y)$
$(j=1, 2, 3)$ of the second kind that have their only pole at $\infty$
such that the fundamental $2$-form of the second kind is given
as\footnote{Since $x$ and $y$ are related, we do not use $\partial$.},
\begin{align}
\Omega((x,y),(z,w))=\left({\frac{d}{dx} \Sigma((z,w),(x,y))+\sum_{k=1}^3
\frac{\omega_k(z,w)}{dz}\frac{\eta_k(x,y)}{dx}}\right)dxdz. 
\lbl{eq3.7} 
\end{align}
The set of differentials $\{\eta_1$, $\eta_2$, $\eta_3\}$ is
determined modulo the space spanned by the $\omega_j$s of
{\rm(\ref{eq1.2})}. 
\end{lemma}
\begin{proof}%%%%%%%%%%%%%%%%%%%%%%%%%%%%%%%%%%%%%%%%%%%%%%%%%%%%%
The 2-form
\begin{equation}
\frac{d}{dz}\Sigma\big((x,y),(z,w)\big) dx dz \lbl{eq3.1.6a}
\end{equation}
satisfies the condition on the poles as a function of $(x,y)$, indeed
one can check that (\ref{eq3.1.6a}) has only a second order pole at
$(x,y)=(z,w)$ whenever $(z,w)$ is an ordinary point or a Weierstrass
point; at infinity the expansion (\ref{u3expansion2}) should be used.
However, the form (\ref{eq3.1.6a}) 
%should have the same property as the form of
%the $(z,w)$-variables.  Nevertheless it 
 has unwanted poles at infinity
as a form in the  $(z,w)$-variables. To restore the symmetry given in
(\ref{eq3.1.6}) we complement (\ref{eq3.1.6a}) by the second term to
obtain (\ref{eq3.7}) with polynomials $\eta_j(x,y)$ which should be
found from (\ref{eq3.1.6}). That results in a system of linear
equations for coefficients of $\eta_j(x,y)$ which is always solvable.
As a result, the polynomials $\eta_i(x,y)$ as well as $F((x,y),(z,w))$
are obtained explicitly.
\end{proof}

\begin{remark}
The 1-form 
\[
\Pi_{(z_1,w_1)}^{(z_2,w_2)}(x,y)= \Sigma((x,y),(z_1,w_1))dx
-  \Sigma((x,y),(z_2,w_2))dx  
\]
is the differential of the third kind,  with first order poles at points 
$(x,y)=(z_1,w_1)$ and $(x,y)=(z_2,w_2)$, and residues $+1$ and $-1$ 
correspondingly.  
\end{remark}

\begin{remark}
  The realization of the fundamental 2-form in terms of the
  Schottky-Klein prime-form and $\theta$-functions is given in
  \cite{ba97}, no.272, and the theory based on the $\theta$-functional
  representation is developed in \cite{fa73}. Here we deal with an
  equivalent algebraic representation of the fundamental 2-form which
  goes back to Klein and %show its function.
  exhibit an algebraic expression for it, that is also mentioned by Fay 
  in \cite{fa73} where the prime-form was defined.  
  The above derivation of the fundamental 2-form is done in 
  \cite{ba97}, around pg.\ 194, and it was reconsidered in
  \cite{eel00} for a large family of algebraic curves.  The case of a
  trigonal curve of genus four was developed in \cite{bg06}, pp.~3617--3618.
\end{remark}

It is easily seen that the $\eta_j$ above is written as
\begin{equation}
 \eta_j(x,y)=\frac{h_j(x,y)}{f_y(x,y)}dx,\quad j=1,2,3,
\end{equation}
where $h_j(x,y) \in \mathbb{Q}[\mu_1,\mu_2, \mu_{4},\mu_5,\mu_8,
\mu_3,\mu_6,\mu_9,\mu_{12}][x,y]$, 
and $h_j$ is of homogeneous weight.
%Now we choose the standard  $\eta_j\,$s  uniquely by requiring 
%(see \cite{bg06}, pg.\ 3618)
%\begin{equation}
%  \text{\it the number of terms in  $h_j(x,y)$  is minimal.}  
%\end{equation}

%It is and easily seen that the $\eta_j$ above is written as
%\begin{equation}
% \eta_j(x,y)=\frac{h_j(x,y)}{f_y(x,y)}dx,\quad j=1,2,3,
%\end{equation}
%where $h_j(x,y)\in\bold Q[\mu_1,\mu_2, \mu_{4},\mu_5,\mu_8,
%\mu_3,\mu_6,\mu_9,\mu_{12}][[x,y]]$, and $h_j$ is of homogeneous

The differentials $\eta_j$ are defined modulo the space of 
holomorphic differentials with the same weight, 
but it is possible to choose the standard $\eta_j\,$s uniquely 
by requiring that for each $j=1$, $2$, $3$ the polynomial $h_j(x,y)$ 
do not contain monomials corresponding to non-gaps 
with bigger $j$. Moreover there exist {\it precisely} $2g=6$ monomials 
defining standard differentials, 
for more details see \cite{bl05}, Chapter 4.

In particular, straightforward calculations lead to the following
expressions
\begin{align}
\begin{split}
h_3(x,y)&=-x^2,\lbl{denom-eta}\\
h_2(x,y)&=-2xy+\mu_1x^2,\\
h_1(x,y)&=-(5x^2+(\mu_1\mu_2-3\mu_3)x+\mu_2\mu_4+\mu_6)y
  +\mu_{2}y^2+3\mu_1x^3 \\
  & \qquad
  -({\mu_2}^2+2\mu_3\mu_1-2\mu_4)x^2-(\mu_5\mu_2+\mu_6\mu_1+\mu_3\mu_4)x+\tfrac34 \mu_1 f_x(x,y)\\ & \qquad
      -\left(\tfrac13\mu_2-\tfrac14{\mu_1}^2\right)f_y.
\end{split}
\end{align}
The orders of monomials defining standard differentials are printed in
bold:
\[ \overline{\mathbf{0}},\;1,\; 2,\;\overline{ \mathbf{3},\;
  \mathbf{4}},\; 5,\; \overline{\mathbf{6},\;\mathbf{ 7},\; 8,\;
  9,\;\mathbf{10},\ldots},
\]
these can be written as $3i+4j$, $0\leq i\leq 2$, $0\leq j\leq 1 $.
We remark that the last two terms in the definition of $h_1(x,y)$
are chosen to provide the standard differentials described above.
The polynomial $F\big((x,y),(z,w)\big)$ in (\ref{realization}) is of
homogeneous weight (weight $-8$), and is given explicitly in Appendix
A.

%%%%%%%%%%%%%%%%%%%%%%%%%%%%%%%%%%%%%%%%%%%%%%%%%%%%%%%%%%%%%%%%%%%%%%%%%%%%%%%
\section{Purely trigonal curve of degree four} %%%%%% Section 6 %%%%%%%%%%%%%%%
%%%%%%%%%%%%%%%%%%%%%%%%%%%%%%%%%%%%%%%%%%%%%%%%%%%%%%%%%%%%%%%%%%%%%%%%%%%%%%%
\label{PureTrig}
In Section \ref{Add2} of this paper, we restrict ourselves to the curve
\begin{equation}
   C:\, y^3  = x^4 + \mu_3 x^3 + \mu_6 x^2 + \mu_9 x + \mu_{12}
   \lbl{eq5.1}
\end{equation}
specialized from (\ref{eq1.1}).  We also restrict results given
in Appendix B to this case to save space.  This curve is called a {\it
  purely trigonal curve} of degree four.  Equivalently we can
represent the curve (\ref{eq1.1}) in the form
\begin{equation}
   C:\, y^3  = \prod_{k=1}^4(x-a_k),
   \lbl{eq5.1a}
\end{equation}
and evaluate the discriminant $D$ according to (\ref{discriminant}) as
\begin{equation} 
  D=\prod_{1\leq i<j \leq 4 } (a_i-a_j)^4.
\end{equation}
% (The sign $\pm1$ is NOT unique!)
The curve $C$ is smooth if and only if $a_i\neq a_j$ for all $i,j=1,
\ldots 4$.  While we assume this to be case, results in the singular
cases are obtained by suitable limiting process.

For the curve (\ref{eq5.1}), the basis (\ref{eq1.2}) of differential
forms of first kind and the function $\Sigma$ in (\ref{eq2.3}) can be
written explicitly as
\begin{equation}
   \omega_1=\frac{dx}{3y^2},  \quad
   \omega_2=\frac{xdx}{3y^2},  \quad
   \omega_3=\frac{ydx}{3y^2}=\frac{dx}{3y},
\end{equation}
and
\begin{equation}
   \Sigma\big((x,y),(z,w)\big)
   =\frac{y^2+yw+w^2}{3(x-z)y^2},
\end{equation} 
respectively. 
The function   $\sigma(u)$  is defined by using these.   
Let 
\[
\zeta=e^{2\pi \sqrt{-1}/3}.
\]
The curve $C$ has an
automorphism $(x,y)\mapsto (x, {\zeta}y)$, and for
$u=(u_1,u_2,u_3)\in\kappa^{-1}\iota(C)$, ${\zeta}^j$ acts by
\begin{equation}
   [{\zeta}^j]u=({\zeta}^ju_1,{\zeta}^ju_2,{\zeta}^{2j}u_3)
                         =\int_{\infty}^{(x,\,{\zeta}^jy)}(du_1,du_2,du_3).
   \lbl{eq5.2}
\end{equation}
This action naturally induces an action on $\kappa^{-1} \Theta^{[k]},
(k=2,3,\dots)$, implying that the set $\Theta^{[k]}$ is
stable under the action of $[{\zeta}^j]$.

%%%%%%%%%%%%%%%%%%%%%%%%%%%%%%%%%%%%%%%%%%%%%%%%%%%%%%%%%%%%%%%%%%%%%%%%%
\section{The $\sigma$-function} %%%%%%%%%%%%%% Section 3 %%%%%%%%%%%%%%%%%%%
%%%%%%%%%%%%%%%%%%%%%%%%%%%%%%%%%%%%%%%%%%%%%%%%%%%%%%%%%%%%%%%%%%%%%%%%%
\label{sigma}
We  construct here the {\it $\sigma$-function} 
\begin{equation}
   \sigma(u)=\sigma(u_1,u_2,u_3)
   \lbl{eq2.1}
\end{equation} 
associated with $C$ for $u \in {\mathbb C}^3$ (see also \cite{bel97},
Chap.1).  We choose closed paths
\begin{equation}
   \alpha_i, \beta_j  \ (1\leqq i, j\leqq 3)
   \lbl{eq2.2}
\end{equation} 
on $C$ which generate $H_1(C,\mathbb{Z})$ such that their
intersection numbers are
$\alpha_i\cdot\alpha_j=\beta_i\cdot\beta_j=0$,
$\alpha_i\cdot\beta_j=\delta_{ij}$.  

Define the period matrices by
\begin{equation}
   \left[\,\omega'  \ \omega''  \right]= 
\left[\int_{\alpha_i}\omega_j \quad \int_{\beta_i}\omega_j\right]_{i,j=1,2,3},
   \,\,
   \left[\,\eta'  \ \eta''  \right]= 
\left[\int_{\alpha_i}\eta_j \quad \int_{\beta_i}\eta_j\right]_{i,j=1,2,3}.
   \lbl{eq2.5}
  \end{equation} 
We can combine these two matrices into 
\begin{equation}
   M=\left[\begin{array}{cc}\omega' & \omega'' \\ \eta' & \eta''
     \end{array}\right].
   \lbl{eq2.6}
\end{equation} 
Then  $M$  satisfies 
\begin{equation}
   M\left[\begin{array}{cc} & -1_3 \\ 1_3 & \end{array}\right]{}^t {M}
   =2\pi\sqrt{-1}\left[\begin{array}{cc} & -1_3 \\ 1_3 &
     \end{array}\right].
   \lbl{eq2.7}
\end{equation} 
This is the {\it generalized Legendre relation} (see (1.14) on p.\,11
of \cite{bel97}).  In particular, ${\omega'}^{-1}\omega''$ is a
symmetric matrix.  We know also that
\begin{equation}
   \mathrm{Im}\,({\omega'}^{-1}\omega'') \qquad
   \mbox{is positive definite.} \qquad\qquad\  \lbl{positive-definiteness}
\end{equation}
By looking at (\ref{eq1.2}), we see the canonical divisor class of $C$
is given by $4\infty$ and we are taking $\infty$ as the point, the
Riemann constant is an element of $\big(\frac12\mathbb{Z}\big)^6$ (see
\cite{mu83}, Coroll.3.11, p.166).  Let
\begin{equation}
   \delta:=\left[\begin{array}{cc}\delta'\ \\
       \delta''\end{array}\right]\in \left(\tfrac12\mathbb{Z}\right)^{6}
   \lbl{eq2.9} %3.15
\end{equation} 
be the theta characteristic which gives the Riemann constant with respect to 
the base point $\infty$ and to the period matrix $[\omega'\ \omega'']$. 
Note that we use $\delta',\delta''$ as well as $n$ in (\ref{def-sigma})
as columns, to keep the notation a bit simpler.
We define
\begin{equation}
\begin{aligned}
   \sigma(u)&=\sigma(u;M)=\sigma(u_1,u_2,u_3;M) \\
   &=c\,\text{exp}(-\tfrac{1}{2}u\eta'{\omega'}^{-1}\ ^t\negthinspace u)
   \vartheta\negthinspace
   \left[\delta\right]({\omega'}^{-1}\ ^t\negthinspace u;\ 
{\omega'}^{-1}\omega'') \lbl{def-sigma}\\
   &=c\,\text{exp}(-\tfrac{1}{2}u\eta'{\omega'}^{-1}\ ^t\negthinspace u)
  \times\\ & \qquad \times \sum_{n \in \mathbb{Z}^3} \exp \big[2\pi i\big\{
   \tfrac12 \ ^t\negthinspace (n+\delta'){\omega'}^{-1}\omega''(n+\delta')
   + \ ^t\negthinspace (n+\delta')({\omega'}^{-1}\,^tu+\delta'')\big\}\big], 
\end{aligned}
\end{equation}
where 
\begin{equation}
    c=\frac{1}{\sqrt[8]{D}}\bigg(\frac{\pi^3}{|\omega'|}\bigg)^{1/2}
   \lbl{sigma-const}
\end{equation}
with $D$ from (\ref{discriminant}).  
Here the choice of a root of (\ref{sigma-const}) is explained 
in the remark \ref{varepsilon} below.  
The series (\ref{def-sigma}) converges 
because of (\ref{positive-definiteness}).

In what follows, for a given $u\in\mathbb{C}^3$, we denote by $u'$ and
$u''$ the unique elements in $\mathbb{R}^3$ such that
\begin{equation}
   u=u'\omega'+u''\omega''.
   \lbl{eq2.12}
\end{equation}
Then for $u$, $v\in\mathbb{C}^3$, and $\ell$
($=\ell'\omega'+\ell''\omega''$) $\in\Lambda$, we define
\begin{align}
  L(u,v)    &:={u}(\eta'{}^tv'+\eta''{}^tv''),\nonumber \\
  \chi(\ell)&:=\exp[\pi\sqrt{-1}\big(2({\ell'}\delta''-
  {\ell''}\delta') +{\ell'}{}^t\ell''\big)] \ (\in \{1,\,-1\}).
   \lbl{eq2.13}
\end{align} 
In this situation the most important properties of $\sigma(u;M)$
are as follows: 
\begin{lemma}%Lemma 2.14}
  The function $\sigma(u)$ is an entire function.  For all
  $u\in\mathbb{C}^3$, $\ell\in\Lambda$ and
  $\gamma\in\mathrm{Sp}(6,{\mathbb{Z}})$, we have \lbl{L2.14}
\begin{align}
\sigma(u+\ell;M) & =\chi(\ell)\sigma(u;M)\exp
L(u+\tfrac12\ell,\ell),\lbl{L2.14.1}\\
\sigma(u;\gamma
M)&=\sigma(u;M),\lbl{L2.14.2}\\
  u\mapsto\sigma(u;M)& \text{ has zeroes of
order $1$ along } \Theta^{[2]}, \lbl{L2.14.3}\\
\sigma(u;M)&=0 \iff u\in\Theta^{[2]}.\lbl{L2.14.4}
\end{align}
\end{lemma}

\begin{proof}% Proof}.
  The function $\sigma$ is clearly entire from its definition and from
  the known property of theta series.  The formula (\ref{L2.14.1}) is
  a special case of the equation from \cite{ba97} (p.286 in the 1995
  reprint, $\ell$.22).  The statement (\ref{L2.14.2}) is easily shown
  by using the definition of $\sigma(u)$ since $\gamma$ corresponds to
  changing the choice of the paths of integration given in
  (\ref{eq2.5}).  The statements (\ref{L2.14.3}) and (\ref{L2.14.4})
  are explained in \cite{ba97}, (p.252).  These facts are partially
  described also in \cite{bel97}, (p.12, Th.1.1 and p.15).
\end{proof}

\begin{lemma}% 2.15
  The function $\sigma(u)$ is either odd or even, i.e. \lbl{parity}
  \begin{equation}
  \sigma([-1]u)=-\sigma(u)\quad \mathrm{or}\quad 
  \sigma([-1]u)= \sigma(u). \lbl{L2.14.0}
  \end{equation}
\end{lemma}

\begin{proof}
  We fix a matrix $M$ satisfying (\ref{eq2.7}) and
  (\ref{positive-definiteness}).  Therefore the bilinear form $L(\ ,\
  )$ is fixed.  Then the space of the solutions of (\ref{L2.14.1}) is
  one dimensional over $\mathbb{C}$, because the Pfaffian of the
  Riemann form attached to $L(\ ,\ )$ is $1$ (see \cite{on98}, Lemma
  3.1.2 and \cite{la82}, p.93, Th.3.1).  Hence, such non-trivial
  solutions automatically satisfy (\ref{L2.14.2}) and (\ref{L2.14.4});
  while (\ref{L2.14.3}) requires the constant factor to be the same,
  this is guaranteed by the definition of $\sigma$ and the fact that
  (\ref{sigma-const}) is independent of $\gamma$.  In this sense,
  (\ref{L2.14.1}) characterizes the function $\sigma(u)$ up to a
  constant, which depends only on the $\mu_j$s.  Now considering the
  loop integrals for $\omega$ in the reverse direction, we see that
  \[
  [-1]\Lambda=\Lambda.
  \]
Hence $u\mapsto \sigma([-1]u)$  satisfies (\ref{L2.14.1}) also. 
So there exists a constant $K$ such that 
  \[
  \sigma([-1]u)=K\,\sigma(u).
  \]
Since $[-1]^2$ is trivial, it must be  $K^2=1$.  \lbl{R2.15}
\end{proof}

\begin{remark}
  In fact $\sigma(u)$ is an {\it odd function} as we see in the
  Theorem \ref{L6.1}.
\end{remark}

We need the power series expansion of $\sigma(u)$ with respect to
$u_1$, $u_2$, $u_3$.  To get the expansion, first of all, we need to
investigate Abelian functions given by logarithmic (higher)
derivatives of $\sigma(u)$.  We shall examine this in the next
Section.

%%%%%%%%%%%%%%%%%%%%%%%%%%%%%%%%%%%%%%%%%%%%%%%%%%%%%%%%%%%%%%%%%%%%%%%%%%%%%%
\section{Standard Abelian functions} %%%%%%%% Section 4 %%%%%%%%%%%%%%%%%%%%%%
%%%%%%%%%%%%%%%%%%%%%%%%%%%%%%%%%%%%%%%%%%%%%%%%%%%%%%%%%%%%%%%%%%%%%%%%%%%%%%
\label{Abelian}
\begin{definition}
A meromorphic function  $u\mapsto\mathfrak{P}(u)$  on  $\mathbb{C}^3$  
is called a {\it standard Abelian function} if it is holomorphic 
outside  $\kappa^{-1}(\Theta^{[2]})$  and is multi-periodic, namely, 
if it satisfies
\begin{equation} 
   \mathfrak{P}(u+\omega'n+\omega''m)=\mathfrak{P}(u) 
\end{equation}
for all integer vectors $n,m\in \mathbb{Z}$ and all $u\not\in
\kappa^{-1}(\Theta^{[2]})$.    \lbl{D4.1}
\end{definition}

To realize the standard Abelian functions in terms of the $\sigma$-function, 
we first let 
\begin{equation}
  \Delta_i=\tfrac{\partial}{\partial u_i}-\tfrac{\partial}{\partial v_i}
\end{equation}
for $u=(u_1,u_2,u_3)$ and $v=(v_1,v_2,v_3)$.  This operator occurs in
what is now known as Hirota's bilinear operator, but in fact was
introduced much earlier in the PDE case by Baker (\cite{ba02}, p.151,
\cite{ba07}, p.49) (see also \cite{ee00}).  We define fundamental
Abelian functions on $J$ by
\begin{equation}
 \wp_{ij}(u)=-\tfrac{1}{2\sigma(u)^2}\Delta_i\Delta_j\,\sigma(u)\sigma(v)|_{v=u}
  =-\tfrac{\partial^2}{\partial u_i\partial u_j}\log\sigma(u).
\lbl{wpij}
\end{equation}
It follows from (\ref{L2.14.0}) that these functions are even.  For
the benefit of the reader familiar with the genus one case, we should
point out that the Weierstrass function $\wp(u)$ described in eqn.\
(\ref{WP}) would be written as $\wp_{11}(u)$ in this notation.  It is
clear that they belong to $\varGamma(J, 2\Theta^{[2]})$.  Moreover, we
define
\begin{equation}
   \wp_{ijk}(u)=\tfrac{\partial}{\partial u_k}\wp_{ij}(u),\quad
   \wp_{ijk\ell}(u)=\tfrac{\partial}{\partial u_{\ell}}\wp_{ijk}(u).
   \lbl{eq3.1}
\end{equation}
The three index $\wp$-functions are odd and four index $\wp$ are even.
The functions (\ref{wpij}) and (\ref{eq3.1}) are standard Abelian
functions from Lemma \ref{L2.14}.  Following (and generalizing) Baker
(\cite{ba02}, pg 151, \cite{ba07}, pp.49--50) (see also \cite{bel00},
pp.18--19, or \cite{cn06}), we define
\begin{equation}
  \begin{aligned}
    Q_{ijk\ell}(u)&=-\tfrac{1}{2\sigma(u)^2}\Delta_i\Delta_j\Delta_k\Delta_{\ell}
    \,\sigma(u)\sigma(v)|_{v=u}\\
    &= \wp_{ijk\ell}(u)-2(\wp_{ij}\wp_{k\ell}+\wp_{ik}\wp_{j\ell}+
    \wp_{i\ell}\wp_{jk})(u),
  \end{aligned}
\lbl{eq3.2a}
\end{equation}
which specializes to 
\begin{align*}
&  Q_{ijkk} = \wp_{ijkk} - 2 \wp_{ij}\wp_{kk}-4\wp_{ik}\wp_{jk},\quad
&& Q_{iikk} = \wp_{iikk} - 2 \wp_{ii}\wp_{kk}-4{\wp_{ik}}^2,\\
&  Q_{ikkk} = \wp_{ikkk} - 6\wp_{ik}\wp_{kk},\quad
&& Q_{kkkk} = \wp_{kkkk} - 6{\wp_{kk}}^2.
\end{align*}
A short calculation shows that $Q_{ijk\ell}$ belongs in
$\varGamma(J,\mathcal{O}(2\Theta^{[2]}))$, whereas
$\wp_{ijk\ell}$ belongs in $\varGamma(J,\mathcal{O}(4\Theta^{[2]}))$.
In particular $Q_{1333}$ plays a key role in what follows.

Note that although the subscripts in $\wp_{ijk\ell}$ {\em do} denote
differentiation, the subscripts in $Q_{ijk\ell}$ do {\em not} denote
direct differentiation, and the latter notation is introduced for
convenience only.  This is important to bear in mind when we use
cross-differentiation, for example the $\wp_{ijk\ell}$ satisfy
\[
    \tfrac{\partial}{\partial u_m}\wp_{ijk\ell}(u) 
  = \tfrac{\partial}{\partial u_\ell}\wp_{ijkm}(u),
\]
whereas the $Q_{ijk\ell}$ do not.  The following useful formula
(\ref{klein}) involving fundamental Kleinian functions, for the case
of the general curve (\ref{eq1.1}), was derived in \cite{bel00}.  It
would be helpful for the reader to read \cite{ba98}, p.\ 377, for the
case of hyperelliptic curves. The formula (\ref{klein}) below is
proved similarly.
%%%%%%%%%%%%%%%%%%%%%%%%%%%%%%%%
\begin{prop} Let $u\in\mathbb{C}^3$ and $(x_1,y_1)$, $(x_2,y_2)$, $(x_3,y_3)$  
be Abelian preimages of  $u$, i.e. 
 \[
 u=\int_{\infty}^{(x_1,y_1)}\omega+\int_{\infty}^{(x_2,y_2)}\omega+ 
 \int_{\infty}^{(x_3,y_3)}\omega
 \]
with appropriate paths of the integrals.  Let $(x,y)$ be an arbitrary
point on the curve $C$.  Then, for each $k=1$, $2$, $3$, the following
formula holds \lbl{Tklein}
\begin{equation}
[1 \ x \ y]\Bigg[\wp_{ij}\bigg(\int_{\infty}^{(x,y)}\omega -u \bigg)\Bigg]
\Bigg[\begin{matrix} \ 1 \ \\ x_k \\ y_k \end{matrix}\Bigg]=
 \frac{F\big((x,y),(x_k,y_k)\big)}{(x-x_k)^2}\lbl{klein},
\end{equation}
where  $F\big((x,y),(z,w)\big)$ is a polynomial defined by {\rm(\ref{realization})} 
or {\rm(\ref{polar})}. 
\end{prop} 

\begin{proof}
Using (\ref{L2.14.3}) and relations of differntials of 
the second kind on $C$ with ones of the third kind 
(see \cite{ba97}, p.22,$\ell$.15 and p.22,$\ell$.11), 
we have an equation connecting the theta series appeared in (\ref{def-sigma}) 
and differentials of the third kind (see \cite{ba97}, p.275, $\ell$.$-11$, for example). 
Then such the equation is modified to a form suitable with $\sigma(u)$ 
and the 2-form $\Omega((x,y),(z,w))$  of (\ref{eq3.7}).  
Finally, after taking logarithm of the modified one, 
operating $\frac{\partial^2}{\partial u_i\partial u_j}$ to it gives the desired equation. 
\end{proof}

\begin{prop}
Suppose the $(x_i,y_i)$s and $u$ are related as in Proposition {\rm \ref{Tklein}}.  
Let $(x,y)$ be any one of $(x_i,y_i)$s. 
Then we have infinitely many relations, of homogeneous weight, linear in 
\lbl{inversion-eq}
  \[
  \wp_{ij}(u),\ \wp_{ijk}(u),\ \wp_{ijk\ell}(u), \cdots  \qquad(i,\ j,\ k = 1,\ 2, 3),
  \]
and whose coefficients are polynomials of  $x$, $y$ and $\mu_j$s. 
We list the first three of them of lower weights as follows\,{\rm\upshape:}
%In view of its central importance, we give here the details of first
%few steps in this theory.  The first three terms in the expansion of
%(\ref{klein}) are the following
  \begin{align}
  &\wp_{33}(u)y+\wp_{23}(u)x +\wp_{13}(u)=x^2, \lbl{kl1} \\
  &\left(\wp_{23}(u)+\tfrac13\mu_1\wp_{33}(u)-\wp_{333}(u)\right)y 
   + \big(\wp_{22}(u)-\wp_{233}(u) \notag\\
  &\qquad\qquad + \tfrac13 \mu_1\wp_{23}(u)\big)x
   +\tfrac13\mu_1\wp_{13}(u) +\wp_{12}(u)-\wp_{133}(u) 
  = 2xy-\tfrac23 \mu_1 x^2, \lbl{kl2}\\
  & -3 y^2+\left(\tfrac13 \wp_{33} \mu_2 +\tfrac12 \wp_{3333}
   -\tfrac12\mu_1 \wp_{333}+\tfrac19 {\mu_1}^2\wp_{33} + 2 \mu_1 x-\tfrac32
  \wp_{233}+2 \mu_4 \right) y \notag\\ 
  & \qquad\qquad +\left(\tfrac23 \mu_2 -\tfrac19 {\mu_1}^2\right)
  x^2+(-\tfrac12 \mu_1\wp_{233} + \mu_5 +\tfrac12 \wp_{2333}+\tfrac13
  \wp_{23} \mu_2 +\tfrac19 {\mu_1}^2\wp_{23} - \tfrac32 \wp_{223})x \notag\\ 
  & \qquad\qquad\qquad + \tfrac12 \wp_{1333} + \tfrac13 \mu_2\wp_{13}
    +\mu_8 +\tfrac19{\mu_1}^2\wp_{13} - \tfrac32 \wp_{123}-\tfrac12 \mu_1 \wp_{133}=0.
  \lbl{kl3}
  \end{align}
More equations of this type are available in {\rm\cite{Weier34}}. 
\end{prop}

\begin{proof}
These relations are derived from (\ref{klein}) 
by expanding (\ref{klein}), with respect to a local parameter  $t=x^{-1/3}$, 
in the vicinity of the point at infinity, 
and comparing the principal parts of the poles on both sides
of the relation (\ref{klein}), we find the solution of the Jacobi
inversion problem.   
%The formula (\ref{klein}) generates various differential relations
%between $\wp$-functions and is the key to all that follows. In
%particular, expanding (\ref{klein}), with respect to a local parameter
%of the variable point $t=x^{-1/3}$, in the vicinity of the point at
%infinity, and comparing the principal parts of the poles on both sides
%of the relation (\ref{klein}), we find the solution of the Jacobi
%inversion problem.   In the case of trigonal curves, this was first
%given explicitly for a particular case of the curve (\ref{eq1.1}) in
%\cite{eel00}. We are lead also to a number of equations which can be
%developed into full sets of equations for the $\wp$ functions.  
\end{proof}

\begin{remark}
(1) In the case of trigonal curves, formula of this type was first
  given explicitly for a particular case of the curve (\ref{eq1.1}) in
  \cite{eel00}.
\newline\noindent
(2) We use in the proof of Lemma \ref{L4index} below the first seven
relations in \ref{inversion-eq}.  Namely, those of weight from $-6$
to $-12$.
\end{remark}

The first two of relations in \ref{inversion-eq} give solution of the
Jacobi inversion problem (see also \cite{bel00}):
\begin{cor}
  Suppose the $(x_i,y_i)$s and $u$ are related as in Theorem {\rm
    \ref{Tklein}}.  The solution of the Jacobi inversion problem is
  given by $(x_1,y_1),(x_2,y_2)$,and $(x_3,y_3)$, where these points
  are the set of zeros of the equations {\rm\upshape(\ref{kl1})},
  {\rm\upshape(\ref{kl2})} for $(x,y)$.  \lbl{jacobi-inversion}
\end{cor}
%%%%%%%%%%%%%%%%%%%
We remark that the right hand side of equations (\ref{kl1}),
(\ref{kl2}) are related to the polynomials $h_3(z,w)$ and $h_2(z,w)$
defining the canonical meromorphic differentials $\eta_3(z,w)$ and
$\eta_2(z,w)$.  {\em Further the first equation in {\rm\upshape (4.8)}
  is directly related to the determinant of the matrices constructed
  in {\rm\cite{on05}}, using the algebraic approach developed in
  \,{\rm\cite{mp05}}.}

If we take the resultant of (\ref{kl1}), (\ref{kl2}) with respect to $y$,
we find a cubic equation in $x$ which can be used to substitute for
$x^3$ in terms of lower powers of $x$
\begin{equation}\begin{split}
x^3 &= \tfrac12 \left(3 \wp_{23}+ \mu_1 \wp_{33}-\wp_{333}\right)
x^2+\tfrac12 \left(\wp_{33} \wp_{22}+2 \wp_{13}+\wp_{23}
  \wp_{333}-\wp_{33} \wp_{233}-{\wp_{23}}^2\right) x\\
& \qquad +\tfrac12 \wp_{33}
\wp_{12}-\tfrac12 \wp_{33} \wp_{133}-\tfrac12 \wp_{13}
\wp_{23}+\tfrac12 \wp_{13} \wp_{333}. \lbl{z3}
\end{split}
\end{equation}
If we now take the resultant of (\ref{kl1}), (\ref{kl3}) with respect
to $y$, we get a quartic in $x$ which can be reduced to a quadratic by
repeated use of (\ref{z3}).  This quadratic in $x$ is not further
reducible A quadratic equation in $x$ has at most only two solutions
and $u$ has three free variables.  Hence each the coefficients of $1$,
$x$, $x^2$ of the quadratic must all be identically zero.
Furthermore, each coefficient can be split into two parts which are
even and odd under the reflection (\ref{eq1.6}), and each of these
parts must vanish.  So each term of order higher than two in the
expansion of (\ref{klein}) can give up to six separate equations
involving the $\wp$ functions.  The simplest two arising from the
resultant of (\ref{kl1}), (\ref{kl3}) are
\begin{align}
  \wp_{222} - 2\wp_{33}\wp_{233} + 2\wp_{23}\wp_{333}
  - \mu_2\wp_{233} + \mu_3\wp_{333} + \mu_1\wp_{223}&=0, \lbl{3index1}\\
  \wp_{23}\wp_{233} - 2 \wp_{33}\wp_{223} + \wp_{333}\wp_{22}
  + 2\wp_{133} + \mu_1(\wp_{23}\wp_{333}-\wp_{33} \wp_{233})&=0,
  \lbl{3index2}
\end{align}
where  $\wp_{ij}=\wp_{ij}(u)$  and  $\wp_{ijk}=\wp_{ijk}(u)$.

%%%%%%%%%%%%%%%%%%%%%%%%%%%%%%%%%%%%%%%%%%%%%%%%%%%%%%%%%%%%%%%%%%%%%%%%%%%
\section{Equations satisfied by the Abelian functions for the general
  trigonal case} %%%
%%%%%%%%%%%%%%%%%%  Section 8 %%%%%%%%%%%%%%%%%%%%%%%%%%%%%%%%%%%%%%%%%%%%
\label{PDEs}
We can use the expansion of (\ref{klein}) as described in the
discussion following Theorem \ref{Tklein} to derive various equations
which the Abelian functions defined by (\ref{eq3.1}) and
(\ref{eq3.2a}) must satisfy.  We consider first the 4-index equations,
the generalizations of $\wp''=6\wp^2-\tfrac12g_2$ in the cubic
(genus 1) case.
\begin{lemma}
The $4$-index functions $\wp_{ijk\ell}$ associated with {\rm (\ref{eq5.1})}
satisfy the following relations\,{\rm :}\lbl{L4index}
\begin{align*}
\wp_{3333}& =6{\wp_{33}}^2 +{\mu_1}^2\wp_{33} -3 \wp_{22} 
  + 2 \mu_1\wp_{23} - 4\mu_2 \wp_{33} -2 \mu_4,\\
\wp_{2333}& =6 \wp_{23}\wp_{33} + {\mu_1}^2\wp_{23} 
  + 3\mu_3\wp_{33} - \mu_2\wp_{23} -\mu_5 -\mu_1\wp_{22},\\
\wp_{2233}& = 4{\wp_{23}}^{2}+ 2 \wp_{33}\wp_{22}
  +\mu_1\mu_3\wp_{33} - \mu_2\wp_{22} + 2 \mu_6 
  + 3\mu_3\wp_{23} + \mu_1\mu_2\wp_{23} + 4 \wp_{13},\\
\wp_{2223}& =6 \wp_{22}\wp_{23} + 4 \mu_1\wp_{13} 
  + \mu_1\mu_3\wp_{23} + \mu_2\mu_3\wp_{33} + 2 \mu_3\mu_4 
  + {\mu_2}^2\wp_{23} + 4 \mu_4\wp_{23} + 3 \mu_3\wp_{22}\\
  &\qquad + 2\mu_1\mu_6 + \mu_2\mu_5 - 2 \mu_5\wp_{33},\\
\wp_{2222}& = 6{\wp_{22}}^2 -2 \mu_2\mu_3\wp_{23} 
  + \mu_1\mu_2\mu_5 + 2 \mu_1\mu_3\mu_4 + 24\wp_{13}\wp_{33}
  + 4{\mu_1}^2\wp_{13} - 4\mu_2\wp_{13} - 4 \wp_{1333} \\
  &\qquad + 4\mu_5\wp_{23} +2 {\mu_1}^2\mu_6 - 2 \mu_2\mu_6
  + \mu_3\mu_ 5 - 3{\mu_3}^2\wp_{33} + 12\mu_6\wp_{33}
  + 4 \mu_4\wp_{22} \\
  &\qquad + {\mu_2}^2\wp_{22} +  4\mu_1\mu_3\wp_{22},\\
\wp_{1233}& = 4 \wp_{13}\wp_{23} +2 \wp_{33}\wp_{12} 
  - 2\mu_1\wp_{33}\wp_{13} - \tfrac13{\mu_1}^{3}\wp_{13} 
  + \tfrac13\mu_1\wp_{1333} +\tfrac13{\mu_1}^2\wp_{12}  
  + 3 \mu_3\wp_{13} \\
  &\qquad +\tfrac13 \mu_1\mu_8 +\tfrac43\mu_1\mu_2\wp_{13}
  - \mu_2\wp_{12} + \mu_9,\end{align*}\begin{align*}
\wp_{1223}& = 4 \wp_{23}\wp_{12} + 2 \wp_{13}\wp_{22} 
  - 2\mu_2\wp_{33}\wp_{13} - 2\mu_8\wp_{33} 
  - \tfrac23\mu_8\mu_2 +\tfrac13\mu_2\wp_{1333} + 3 \mu_3 \wp_{12}  
  + 4\mu_4\wp_{13} \\
  &\qquad + \tfrac43 {\mu_2}^2\wp_{13} - 2 \wp_{11} 
  - \tfrac13{\mu_1}^2\mu_2\wp_{13} +\tfrac13 \mu_1\mu_2\wp_{12} 
  + \mu_1\mu_3\wp_{13},\\
\wp_{1222}& =6 \wp_{22}\wp_{12} + 6\mu_9\wp_{33} - \mu_3\wp_{1333}
  + 4 \mu_5\wp_{13} + {\mu_2}^2\wp_{12} - \mu_2\mu_9 
  + 4 \mu_ 4\wp_{12} - 2 \mu_1\wp_{11}  \\
  &\qquad + 6\mu_3\wp_{33}\wp_{13} -3\mu_2\mu_3\wp_{13}
  + {\mu_1}^2\mu_3\wp_{13} + 3 \mu_1\mu_3\wp_{12} - \mu_1\mu_2\mu_8,\\
\wp_{1133}& =4{\wp_{13}}^2 +2 \wp_{33}\wp_{11} - \mu_9\wp_{23}
  + 2\mu_6\wp_{13} +\mu_8\wp_{22} - \mu_5\wp_{12} 
  + \tfrac23 \mu_4\wp_{1333} + \tfrac23 \mu_4\mu_8 \\
  &\qquad +2\mu_2\mu_8\wp_{33} - 4\mu_4\wp_{13}\wp_{33} 
  + \tfrac23\mu_2\mu_4\wp_{13} + \mu_1\mu_9\wp_{33} - \mu_1\mu_8\wp_{23}
  + \mu_1\mu_5\wp_{13}  \\
  &\qquad -\tfrac23{\mu_1}^2\mu_4\wp_{13} + \tfrac23 \mu_1\mu_4\wp_{12},\\
\wp_{1123}& =4 \wp_{12}\wp_{13} + 2 \wp_{23}\wp_{11}
  + 2\mu_3\mu_4\wp_{13} - \mu_3\mu_8\wp_{33} - 2\mu_5\wp_{13}\wp_{33}
  + \mu_2\mu_8\wp_{23} + \tfrac43\mu_2\mu_5\wp_{13} \\
  &\qquad -\mu_9\wp_{22} + 2\mu_6\wp_{12} + \tfrac13\mu_5\wp_{1333}
  + \tfrac13\mu_5\mu_8 +\mu_1\mu_9\wp_{23} -\tfrac13{\mu_1}^2\mu_5\wp_{13}
  +\tfrac13 \mu_1\mu_5\wp_{12},\\
\wp_{1122}& = 4{\wp_{12}}^2 +2 \wp_{11}\wp_{22} +\tfrac23
  {\mu_1}^2\mu_6\wp_{13} + \tfrac43 \mu_1\mu_6\wp_{12} 
  + \mu_3\mu_9\wp_{33} + \mu_2\mu_9\wp_{23} 
  + 8\mu_{12}\wp_{33}  \\
  &\qquad + 2 \mu_3\mu_4\wp_{12}-\tfrac23\mu_6\wp_{1333} 
  + 4\mu_8\wp_{13} -\tfrac23\mu_6\mu_8 + 4\mu_6\wp_{33}\wp_{13} 
  - \mu_3\mu_8\wp_{23} + \mu_3\mu_5\wp_{13} \\
  &\qquad - \tfrac83 \mu_2\mu_6\wp_{13} + \mu_2\mu_8\wp_{22} 
  + \mu_2\mu_5\wp_{12},\\
\wp_{1113}& = 6 \wp_{13}\wp_{11} + 6 \mu_2\mu_8\wp_{13}
  - 2\mu_2\mu_{12}\wp_{33} - {\mu_1}^2\mu_8\wp_{13} 
  + 4\mu_1\mu_{12}\wp_{23} + \mu_1\mu_8\wp_{12} + \mu_5\mu_9\wp_{33} \\
  &\qquad + {\mu_5}^2\wp_{13} - 2 \mu_4\mu_9\wp_{23} + \mu_1\mu_9\wp_{13} 
  - 6\mu_8\wp_{33}\wp_{13} - 2\mu_6\mu_8\wp_{33} + \mu_8\wp_{1333} 
  - 4\mu_4\mu_{12} \\
  &\qquad + 3 \mu_9\wp_{12} - 6 \mu_{12}\wp_{22} - \mu_5\mu_8\wp_{23}
  + 4 \mu_4\mu_6\wp_{13},\\
\wp_{1112}& = 6 \wp_{12}\wp_{11} + 6\mu_3\mu_{12}\wp_{33}
  + 3 \mu_3\mu_8\wp_{13} - 2\mu_6\mu_8\wp_{23} - \mu_1{\mu_8}^2 
  + 5\mu_2\mu_8\wp_{12} + 4\mu_2\mu_{12}\wp_{23}  \\
  &\qquad - 2 \mu_1\mu_{12}\wp_{22} + 4 \mu_4\mu_6\wp_{12}
  - \mu_5\mu_8\wp_{22} + {\mu_5}^2\wp_{12} + 4\mu_5\mu_{12}
  - \mu_9\wp_{1333} - 4 \mu_1\mu_{12}\mu_4 \\
  &\qquad + {\mu_ 1}^2\mu_9\wp_{13} + 3 \mu_1\mu_9\wp_{12} 
  - 2 \mu_4\mu_9\wp_{22} +\mu_5\mu_9\wp_{23} -4 \mu_2\mu_9\wp_{13} 
  + 6 \mu_9\wp_{13}\wp_{33} - 3 \mu_8\mu_9,\\
\wp_{1111}& =6{\wp_{11}}^2 + 4 \mu_4\mu_9\wp_{12} - 8{\mu_4}^2\mu_{12}
  - 2{\mu_2}^2\mu_4\mu_{12} -3 {\mu_8}^2\wp_{22} - 2 \mu_4{\mu_8}^2
  + {\mu_5}^2\wp_{11} -3{\mu_9}^2\wp_{33} \\
  &\qquad - 4\mu_{12}\wp_{1333} + 24 \mu_{12}\wp_{33}\wp_{13} 
  + 12 \mu_5\mu_{12}\wp_{23} + \mu_2 \mu_4\mu_5\mu_9
  - 6 \mu_1\mu_3\mu_4\mu_{12} \\
  &\qquad +\mu_1\mu_2\mu_5\mu_{12} +2 {\mu_6}^2\mu_8 +2{\mu_2}^2{\mu_8}^2
  - \mu_5\mu_6\mu_9 - 2 \mu_5\mu_9\wp_{13} 
  + 4\mu_4\mu_6\wp_{11} + 4\mu_6\mu_8\wp_{13} \\
  &\qquad + 8\mu_2\mu_8\wp_{11} - 6 \mu_2\mu_6\mu_{12}
  - 12 \mu_2\mu_{12}\wp_{13} + 4{\mu_1}^2\mu_{12}\wp_{13} 
  + 2 {\mu_1}^2\mu_6\mu_{12} + 2\mu _8\mu_5\wp_{12} \\
  &\qquad - 6 \mu_8\mu_9\wp_{23} - 12\mu_4\mu_{12}\wp_{22} + \mu_2 {\mu_5}^2\mu_8
  + 2\mu_1\mu_4\mu_6\mu_9 +\mu_1\mu_5\mu_6\mu_8 + 12\mu_6\mu_{12}\wp_{33} \\
  &\qquad + 4\mu_1\mu_9\wp_{11} + 2\mu_3 {\mu_4}^2\mu_9
  + 9\mu_3\mu_5\mu_{12} - 2\mu_1\mu_3{\mu_8}^2 - 6\mu_3\mu_8\mu_9 
  + 2\mu_1\mu_2\mu_8\mu_9 \\
  &\qquad + \mu_3\mu_4\mu_5\mu_8  + 2\mu_2\mu_4\mu_6 \mu_8
  + 2\mu_2{\mu_9}^2. 
\end{align*}
\end{lemma}
\begin{proof}
  Many of these relations follow from the sets of equations generated
  from the first seven terms of the expansion of (\ref{klein}) as
  indicated in Proposition \ref{inversion-eq} by a similar argument as
  that explained at the end of the previous Section.  Others can be
  derived making use of derivatives of the equations in Lemma
  \ref{L3index}, or products of these equations with three index
  expressions $\wp_{ijk}$, working in a self-consistent way from
  higher to lower weights.  The calculations are somewhat long and
  tedious and much facilitated by heavy use of Maple.  Full Maple
  worksheets are available on request from the authors
\end{proof}

\begin{remark} The complete set of the four-index relations for
  $\wp$-functions for genus three was derived by Baker \cite{ba02} in
  the hyperelliptic case only. As far as we know, the above relations
  are new, and a comparison with Baker's relations is of interest.
\end{remark}

\begin{remark}
  With the use of (\ref{eq3.2a}), these equations can be written in a
  slightly more compact form involving the $Q_{ijk\ell}$ functions.
  For example, the sixth equation (for $\wp_{2222}$) becomes
\begin{align*}
  Q_{2222}& = -2\mu_2\mu_3\wp_{23} + \mu_1\mu_2\mu_5 
  + 2\mu_1\mu_3\mu_4 + 4{\mu_1}^2\wp_{13} 
  - 4\mu_2\wp_{13} - 4 Q_{1333} + 4\mu _5\wp_{23} \\
  &\qquad +2 {\mu_1}^2\mu_6 - 2\mu_2\mu_6
  + \mu_3\mu_ 5 - 3{\mu_3}^2\wp_{33} + 12\mu_6\wp_{33}
  + 4\mu_4\wp_{22} + {\mu_2}^2\wp_{22} + 4\mu_1\mu_3\wp_{22}.
\end{align*}
The importance of this switch to the $Q$ variables is that the
equations become {\em linear} in the $Q_{ijk\ell}$ and the 2-index
$\wp_{ij}$.  An alternative way of looking at this is that the
equations in Lemma \ref{L4index} have only second-order poles in
$\sigma$.
\end{remark}
\begin{remark}
  The first relation in Lemma \ref{L4index}, after differentiating
  twice with respect to $u_3$, becomes the Boussinesq equation for the
  function $\wp_{33}$ (see \cite{bel00,eel00}).
\end{remark}
\begin{lemma}
The $3$-index functions $\wp_{ijk}$ associated with {\rm (\ref{eq5.1})}
satisfy a number of bi-linear relations {\rm (}linear in both $3$-index and
$2$-index functions{\rm )}.  These have no analogue in the genus $1$ case.
For example, in decreasing weight, starting at $-6$ we have \lbl{L3index}
\begin{alignat*}{2}
  & {-}2\,\wp_{33}\wp_{233}+2\,\wp_{23}\wp_{333}+\mu_3\wp_{333}
       +\mu_2\wp_{233}-\mu_1\wp_{223}+\wp_{222}=0,\quad                   &[-6]\\
  & {-}2\,\wp_{33}\wp_{223}+\mu_1\wp_{33}\wp_{233}+\wp_{23}\wp_{233}
       -\mu_1\wp_{23}\wp_{333}+\wp_{333}\wp_{22}+2\,\wp_{133}=0,        \quad &[-7]\\
  & {-}2\,\wp_{23}\wp_{223}+4\,\wp_{22}\wp_{233}+4\,\mu_1\wp_{133}
       +\mu_3\mu_2\wp_{333}+{\mu_2}^2\wp_{233}                      &\\
       & \qquad +4\,\mu_4\wp_{233}-2\,\mu_5\wp_{333}
       -2\,\wp_{33}\wp_{222}+\mu_2\wp_{222}
       -\mu_1\mu_2\wp_{223}-4\,\wp_{123}=0,                         \quad &[-8]\\
  & 3\,\mu_1\mu_3\wp_{223}-3\,\mu_2\mu_3\wp_{233}
       -24\,\wp_{33}\wp_{133}+24\,\wp_{13}\wp_{333}+12\,\wp_{122}
       -12\,\mu_1\wp_{123}                                                    &\\
       & \qquad +12\,\mu_2\wp_{133}-3\,\mu_3\wp_{222}
       +6\,\mu_5\wp_{233}-3\,{\mu_3}^2\wp_{333}
       +12\,\mu_6\wp_{333}-6\,\wp_{23}\wp_{222}+6\,\wp_{22}\wp_{223}=0, \quad &[-9]\\
  & 2\,\wp_{33}\wp_{123}-\mu_1\wp_{33}\wp_{133}+\mu_1\wp_{13}\wp_{333}
            +\wp_{23}\wp_{133}-\wp_{12}\wp_{333}-2\,\wp_{13}\wp_{233}=0,         \quad &[-10]\\
  & \wp_{113}+\wp_{13}\wp_{223}-2\,\mu_4\wp_{133}+\wp_{33}\wp_{122}
           -\wp_{22}\wp_{133}-\wp_{12}\wp_{233}+\mu_8\wp_{333}
           -\mu_2\wp_{133}\wp_{33}                                            &\\
  & \qquad -\mu_1\wp_{13}\wp_{233}+\mu_2\wp_{13}\wp_{333}
           +\mu_1\wp_{133}\wp_{23}=0,                                   \quad &[-11]\\
  & {-}\wp_{112}-3\,\mu_9\wp_{333}+\wp_{13}\wp_{222}-\wp_{12}\wp_{223}
           -2\,\wp_{22}\wp_{123}-2\,\mu_5\wp_{133}
           +\mu_1\wp_{113}+2\,\wp_{23}\wp_{122}                              &\\
  & \qquad  -\mu_8\wp_{233}-\mu_2\wp_{13}\wp_{233}
           +3\,\mu_3\wp_{33}\wp_{133}-3\,\mu_3\wp_{13}\wp_{333}
           +\mu_2\wp_{23}\wp_{133}=0,                                   \quad &[-12]\\
  & 8\,\mu_4\wp_{133}\wp_{33}-8\,\mu_4\wp_{13}\wp_{333}
           -4\,\mu_2\mu_4\wp_{133}+2\,\mu_1\mu_9\wp_{333}
           -2\,\mu_1\mu_8\wp_{233}+2\,\mu_1\mu_5\wp_{133}         &\\
        & \qquad  +4\,\mu_1\mu_4\wp_{123}+4\,\mu_8\mu_2\wp_{333}
           +3\,\mu_3\wp_{13}\wp_{233}-3\,\mu_3\wp_{23}\wp_{133}
           -\mu_1\wp_{112}+3\,\wp_{12}\wp_{222}                               &\\
        & \qquad +4\,\wp_{11}\wp_{333}-2\,\mu_6\wp_{133}-3\,\wp_{122}\wp_{22}
          -4\,\mu_4\wp_{122}+\mu_9\wp_{233}
          +2\,\mu_8\wp_{223}-8\,\wp_{33}\wp_{113}                             &\\
        & \qquad +4\,\wp_{13}\wp_{133}-2\,{\mu_1}^2\wp_{113}
          +2\,\mu_2\wp_{113}-2\,\mu_5\wp_{123}=0,
          \quad &[-13]
\end{alignat*}
\begin{alignat*}{2}
  & 4\,\wp_{123}\wp_{13}+4\,\mu_4\wp_{23}\wp_{133}
          +\mu_3\mu_8\wp_{333}-2\,\mu_5\wp_{33}\wp_{133}
          +2\,\mu_5\wp_{13}\wp_{333}+\mu_2\mu_8\wp_{233}
          +\mu_8\wp_{222}                                                      &\\
        & \qquad -4\,\wp_{12}\wp_{133}-2\,\wp_{23}\wp_{113}
          + 2\,\wp_{33}\wp_{112}
          - 4\,\mu_4\wp_{13}\wp_{233}-\mu_1\mu_8\wp_{223}=0,       \quad &[-14]\\
  & {-}\mu_9\wp_{222}+\mu_1\mu_9\wp_{223}+4\,\wp_{13}\wp_{122}
           +2\,\wp_{23}\wp_{112}-2\,\wp_{113}\wp_{22}
           -\mu_3\mu_9\wp_{333}-\mu_2\mu_9\wp_{233}                            &\\
        & \qquad +2\,\mu_5\wp_{23}\wp_{133}-8\,\mu_{12}\wp_{333}
           -4\,\mu_8\wp_{133}-4\,\mu_6\wp_{13}\wp_{333}
           +4\,\mu_6\wp_{33}\wp_{133}-4\,\wp_{12}\wp_{123}                     &\\
        & \qquad -2\,\mu_5\wp_{13}\wp_{233}=0.                     \quad &[-15]
\end{alignat*}
where the number in brackets $[\quad ]$ indicates the weight.  
\end{lemma}
\begin{proof}
  We have already given the first two of these equations in the
  discussion following Theorem \ref{Tklein}.  Some of the others
  follow in the same way from the expansion of (\ref{Tklein}).
  Alternatively, some can be calculated directly by expressing the
  equations in Lemma \ref{L4index} in terms of $\wp_{ijk\ell}$ and
  $\wp_{mn}$ functions, then using cross differentiation on suitably
  chosen pairs of equations.  For example the first relation above for
  $\wp_{222}$ can be derived from
\[
 \frac{\partial}{\partial u_2}\wp_{3333} 
-\frac{\partial}{\partial u_3}\wp_{2333}=0.
\]
\end{proof}
\begin{remark}
  For a fixed weight, these relations are not always unique, for
  example at weight $-11$ we also have the relation \lbl{non-unique}
\begin{equation*}
  \wp_{33}\wp_{122}+2\wp_{23}\wp_{123}+3\wp_{113}
  +\mu_2\wp_{13}\wp_{333}-\mu_2\wp_{33}\wp_{133}
  +\mu_8\wp_{333}-2\wp_{12}\wp_{233}-2\mu_4\wp_{133}
  -\wp_{13}\wp_{223}=0
\end{equation*}
These dual relations arise because in some cases the cross
differentiation can be done in two different ways.  In deriving the
results in this section, it is sometimes required to make use of both
bilinear relations at a given weight to provide enough equations to
solve for the unknowns.  A full list of the known bilinear relations
is given at \cite{Weier34}. 
\end{remark}
\begin{lemma}
  The quadratic expressions in the $3$-index functions $\wp_{ijk}$
  associated with {\rm (\ref{eq5.1})} down to weight $-23$ can be
  expressed in terms of {\rm (}at most cubic{\rm )} relations in the
  $\wp_{mn}$ and $\wp_{1333}$.  For example we have the following five
  relations down to weight $-8$\,{\rm:}\lbl{L33index}
\begin{align*}
{\wp_{333}}^2 & = {\wp_{33}}^2{\mu_1}^2 
  + 2\mu_1\wp_{23}\wp_{33} + {\wp_{23}}^2 + 4 \wp_{13}
  - 4\wp_{33}\wp_{22} + 4{\wp_{33}}^3 - 4\mu_2{\wp_{33}}^2
  - 4 \mu_4\wp_{33},\\
\wp_{233}\wp_{333} & = 2 \mu_3{\wp_{33}}^2 + 4 {\wp_{33}}^2\wp_{23}
  - \mu_1\wp_{33}\wp_{22} - 2 \mu_5\wp_{33} 
  - 2\mu_2\wp_{33}\wp_{23} + {\mu_1}^2\wp_{33}\wp_{23}
  - 2\wp_{12}\\
  & \quad -\wp_{22}\wp_{23}+\mu_1{\wp_{23}}^2 + 2 \mu_1\wp_{13},\\
\wp_{133}\wp_{333} & =-\tfrac13 \mu_1\wp_{33}\wp_{12}
  + \tfrac13{\mu_1}^2\wp_{33}\wp_{13}-\tfrac43\mu_2\wp_{33}\wp_{13}
  + \tfrac23\wp_{33}\wp_{1333}-\tfrac43\mu_8\wp_{33}+\wp_{23}\wp_{12}\\
  & \quad + \mu_1\wp_{13}\wp_{23} - 2 \wp_{13}\wp_{22},\\
\wp_{223}\wp_{333} & = 2\mu_1\wp_{23}\wp_{22} 
  - 2 \mu_2\wp_{33}\wp_{22} + 2\mu_1\mu_4\wp_{23}
  - \mu_1\mu_5\wp_{33} + 2{\wp_{33}}^2\wp_{22} 
  - 2 \mu_4\wp_{22} + 2 \wp_{33}{\wp_{23}}^2\\
  & \quad + \tfrac43{\mu_1}^2 \wp_{13} - \tfrac43\mu_2\wp_{13}
  - \tfrac43 \mu_1\wp_{12} - \tfrac43 \mu_8 - 2{\wp_{22}}^2 
  + \mu_1\mu_2\wp_{33}\wp_{23} + \tfrac23 \wp_{1333}\\
  & \quad +\wp_{23}\wp_{33}\mu_3 + \mu_1\mu_3{\wp_{33}}^2
  - \mu_2{\wp_{23}}^2 -\mu_5 \wp_{23},\\
{\wp_{233}}^2 & = 4 \wp_{33}{\wp_{23}}^2 + 8 \wp_{13}\wp_{33}
  + 4 \mu_3\wp_{33}\wp_{23} - 2\mu_1\wp_{23}\wp_{22} 
  + \tfrac43{\mu_1}^2\wp_{13} - \tfrac43\mu_2\wp_{13}
  + 4 \mu_6\wp_{33}\\
  & \quad +{\mu_1}^2{\wp_{23}}^2 - \tfrac43 \mu_8 
  + {\wp_{22}}^{2} - \tfrac43 \wp_{1333} - \tfrac43\mu_1\wp_{12}.
\end{align*}
The expressions at lower weight quickly become very lengthy.  For the
purely trigonal case we give a list of the known quadratic
expressions in the $3$-index functions up to weight $-15$ in {\rm
  Appendix B}.  The full list for the general $(3,4)$-curve down to
weight $-23$ is available at {\rm\cite{Weier34}}. 
\end{lemma}

\begin{proof}
  The relations can be found using a combination of three types of
  intermediate relations.  One type is from terms in the expansion of
  (\ref{klein}).  Another is to multiply one of the linear three-index
  $\wp_{ijk}$ relations above by another $\wp_{ijk}$ and substitute
  for previously calculated $\wp_{ijk}\wp_{\ell mn}$ relations of
  higher weight.  Yet another is to take a derivative of one of the
  bilinear three-index $\wp_{ijk}$ relations above and to substitute
  the known linear four-index $\wp_{ijk\ell}$ and previously
  calculated $\wp_{ijk}\wp_{\ell mn}$ relations.  Again, we work in a
  self-consistent way from higher to lower weights.  The strategy for
  all the results in this section is to proceed down one weight at a
  time and to derive {\em all} the three types of relations (4-index
  $\wp_{ijk\ell}$, bilinear 2- and 3-index, and quadratic 3-index) at
  a given weight before moving down to the next.  An extra
  complication is that at certain weights some of the intermediate
  calculations can involve {\em quartic} terms in the $\wp_{mn}$ and
  $\wp_{1333}$.  It is always possible to find enough relations to
  eliminate the quartic term up to weight $-23$.  
\end{proof}
\begin{remark}
  (1) These relations are the generalizations of the familiar relation
  $(\wp')^2 = 4 \wp^3-g_2\wp-g_3$ in the genus 1 theory.
  \\
  (2) For equations of weight below $-23$, we have not been able to
  find cubic expressions for the $\wp_{ijk}\wp_{\ell m n}$ terms.  We
  believe it should be possible to explain this using the results of
  Cho and Nakayashuiki \cite{cn06}, and we are currently investigating
  this possibility.
  \\
  (3) The calculations in this section make no use of the expansion of
  the $\sigma$ function, which is given in the next section.
\end{remark}
%%%%%%%%%%%%%%%%%%%%%%%%%%%%%%%%%%%%%%%%%%%%%%%%%%%%%%%%%%%%%%%%%%%
\section{Expansion of the $\sigma$-function} %%% Section 7 %%%%%%%%%%%
%%%%%%%%%%%%%%%%%%%%%%%%%%%%%%%%%%%%%%%%%%%%%%%%%%%%%%%%%%%%%%%%%%%
\label{sigma_expan}
This Section is devoted to show the coefficients of 
the power series expansion of $\sigma(u)$  is a polynomial in $\mu_j$s.   
%%%%%%%%%%%%%%%%%%%%%%%%%%%%%%%%%%Victor 09/16/06

In the Weierstrass formulation of the theory of elliptic functions,
the $\sigma$-function is defined as the power series expansion in the
Abelian variable $u$ with coefficients depending on the Weierstrass
parameters $g_2,g_3$, and related by certain recursive relations. The
extension of Weierstrass theory to arbitrary algebraic curves was
intensively developed in the 19th century and later, its development
being attached to names such as Baker, Bolza, Brioschi, Burkhardt,
Klein, and Wiltheiss.  Some important modern developments of this
theory are due to Buchstaber and Leykin \cite{bl02,bl05} who
give a construction of linear differential (heat-like) operators that
annihilate the $\sigma$-function for any $(m,n)$-curve. In the
hyperelliptic case the operators are sufficient to find the recursion
defining the whole series expansion. The exact analogue of the
Weierstrass recursive series formula is known only for genus two, see
\cite{bl05}, p.68. %[14]
In other cases the detailed results have not yet been developed,
although the general method is provided in the publications mentioned
above. Here we shall give the few first terms of the power series
expansion, obtained by finding the coefficients of the Taylor series
by using the PDEs given in Lemma \ref{L4index}.
%%%%%%%%%%%%%%%%%%%%%%%%%%%%%%%%%%%%%%%%%%%%%%%%%%%%%%%%%%%%%%%%%%%

\begin{theorem} % 8.1
  The function $\sigma(u)$ associated with the general trigonal curve
  {\rm (\ref{eq1.1})} of genus three has an expansion of the following
  form\,{\rm :}\lbl{L6.1}
  \begin{equation}
  \begin{aligned}
    \sigma(u_1,u_2,u_3)
    &=\varepsilon\cdot\big(C_5(u_1,u_2,u_3)+C_6(u_1,u_2,u_3)
    +C_{7}(u_1,u_2,u_3) +\cdots\big),
  \end{aligned}
\end{equation}
where $\varepsilon$ is a non-zero constant and each $C_{j}$ is a
polynomial composed of sums of monomials in $u_i$ of odd total degree
and of total weight $j$ with polynomial coefficient in $\mu_i$s of
total weight $(5-j)$.  Especially, $\sigma(u)$ is an odd function
{\rm\upshape(}see {\rm\upshape\ref{parity}}{\rm\upshape)}.  The first
few $C_j$s are
\begin{align*}
   C_5 & = u_1-u_3\,{u_2}^2+\tfrac1{20}\,{u_3}^5, \qquad\qquad\quad\quad\;\;
   C_6  = \tfrac1{12}\,\mu_1{u_3}^4u_2
         - \tfrac13\,\mu_1{u_2}^3,\\
   C_7 & = \tfrac {1}{504}\,\big({\mu_1}^{2}-3\,\mu_2 \big){u_3}^7 
         + \tfrac16\,\mu_2{u_3}^3{u_2}^2,\quad
   C_8  = \tfrac {1}{360}\,\big({\mu_1}^3+9\,\mu_3-2\,\mu_1\mu_2\big){u_3}^6u_2
         - \tfrac12\,\mu_3{u_3}^2{u_2}^3,\\
   C_9 & = \tfrac {1}{25920}\,\big({\mu_1}^2-3\,\mu_2\big)^2{u_3}^9
         + \tfrac {1}{120}\,\big(2\,\mu_4-{\mu_2}^2
            +{\mu_1}^2\mu_2+6\,\mu_1\mu_3 \big){u_3}^5{u_2}^2\\
       & \quad -\tfrac1{12}\,\big( 4\,\mu_1\mu_3+4\,\mu_4+ \mu_2^2 
       \big) u_3{u_2}^4  + \tfrac1{12}\,\mu_4{u_3}^4 u_1,\\
C_{10} & = \tfrac{1}{20160}\,\big(8\,\mu_1\mu_4 - 54\,\mu_2\mu_3
          + 3\,\mu_1{\mu_2}^2 + 18\,{\mu_1}^2\mu_3+{\mu_1}^5
          - 12\,\mu_5 - 4\,{\mu_1}^3\mu_2\big){u_3}^8u_2\\
       & \quad + \tfrac {1}{72}\,\big(6\,\mu_2\mu_3 + 2\,\mu_1\mu_4
          + \mu_1{\mu_2}^2 + {\mu_1}^2\mu_3 \big){u_3}^4{u_2}^3 \\
       & \quad -\tfrac {1}{60}\,\big( 4\,{\mu_1}^2\mu_3 + \mu_1{\mu_2}^2
          + 4\,\mu_5 + 4\,\mu_1\mu_4 - 2\,\mu_2\mu_3 \big){u_2}^5
          + \tfrac16\,\mu_5{u_3}^3u_2u_1,\\
C_{11} & = -{\tfrac {1}{6652800}}\, 
          \big(18\,\mu_1\mu_2\mu_3+27\,{\mu_1}^4 \mu_2
          - 72\,\mu_6 -3\,{\mu_1}^6 - 24\,\mu_2 \mu_4
          + 16\,{\mu_1}^2\mu_4  - 24\,\mu_1\mu_5 \\ 
       & \qquad + 27\,{\mu_3}^2 + 85\,{\mu_2}^3 - 4\,{\mu_1}^3\mu_3
          - 82\,{\mu_1}^2{\mu_2}^2 \big){u_3}^{11}
          + \tfrac{1}{5040} \, 
          \big(27\,{\mu_3}^2+{\mu_2}^3-6\,\mu_2\mu_4 \\ 
       & \qquad -18\,\mu_1\mu_2\mu_3 + 8\,{\mu_1}^{3}\mu_3
          - 4\,\mu_1\mu_5 + 6\,{\mu_1}^2\mu_4 
          + 12\,\mu_6 + {\mu_1}^4\mu_2 - 3\,{\mu_1}^2{\mu_2}^2\big)
         {u_3}^7{u_2}^2\\
       & \quad -\tfrac {1}{72}\, 
         \big(9\,{\mu_3}^2-{\mu_2}^3-4\,\mu_2\mu_4
          - 2\,\mu_1\mu_2\mu_3 \big){u_3}^3{u_2}^4 \\
       & \quad + \tfrac {1}{360}\,\big(\mu_1\mu_5
          - 4\,\mu_2\mu_4 + {\mu_1}^2\mu_4+3\,\mu_6 \big)
          {u_3}^6u_1 - \tfrac12\,\mu_6{u_3}^2{u_2}^2u_1,\\
C_{12}& = -\tfrac {1}{1814400}\, 
        \big( 27\,\mu_1{\mu_3}^2 - 243\,{\mu_2}^2\mu_3-{\mu_1}^7 
         + 72\,\mu_1\mu_2\mu_4
         - 31\,{\mu_1}^4\mu_3 - 144\,\mu_2\mu_5 - 16\,{\mu_1}^3\mu_4\\
       &\quad
         + 6\,{\mu_1}^5\mu_2 - 10\,{\mu_1}^3{\mu_2}^2
         + 24\,{\mu_1}^2\mu_5 + 4\,\mu_1{\mu_2}^3
         - 72\,\mu_1\mu_6 + 180\,{\mu_1}^2\mu_2\mu_3 \big){u_3}^{10}u_2\\
       & \quad +{\tfrac {1}{2160}}\, \big( 18\,\mu_3\mu_4
         - 2\,\mu_1{\mu_2}^3 + 27\,\mu_1{\mu_3}^2 - 9\,{\mu_2}^2\mu_3
         + {\mu_1}^3{\mu_2}^2 + {\mu_1}^4\mu_3 
         + 6\,{\mu_1}^2\mu_2\mu_3 \\ 
       & \quad  + 2\,{\mu_1}^3\mu_4 + 12\,\mu_1\mu_6 \big){u_3}^6{u_2}^3 
         - \tfrac1{24}\,\mu_3 \big(3\,\mu_1\mu_3+4\,\mu_4 +{\mu_2}^2\big) 
         {u_3}^2{u_2}^5 \\ 
       & \quad + \tfrac {1}{120}\,
           \big(6\,\mu_3\mu_4+2\,\mu_1\mu_6-\mu_2\mu_5+{\mu_1}^2\mu_5\big)
         {u_3}^5u_2u_1
         - \tfrac16\,\big( 2\,\mu_1\mu_6 + 2\,\mu_3\mu_4+\mu_2\mu_5\big)
          u_3{u_2}^3u_1.
\end{align*}
\end{theorem}

\begin{proof} We divide the proof into four parts.
\vskip 5pt
\noindent{\bf Step 1.}  We have already shown in \ref{parity}, that all
  the terms are of total odd degree or even degree.  We first show
  that the expansion contains a term linear in $u_1$, so the expansion
  must be odd.

  Let $B(D)$ be the Brill-Noether matrix for an effective divisor $D$
  of $C$.  Then it is well known that (see for example \cite{on98} or
  \cite{onPEMS})
  \[
  \dim\varGamma(C,\mathcal{O}(D))=\deg D + 1 -\mathrm{rank} B(D),
  \]
  where $\varGamma(C,\mathcal{O}(D))$ is the space of functions on $C$
  whose divisor are larger than or equal to $-D$.  Moreover, for two
  points $P_1$, $P_2$ on $C$,
  $\dim\varGamma(C,\mathcal{O}(P_1+P_2))>1$ if and only if the point
  $\iota(P_1,P_2)\in \Theta^{[2]}$ is a non-singular point of
  $\Theta^{[2]}$ (note that $C$ is of genus $3$).  By checking the
  Brill-Noether matrix $B(P_1+P_2)$, we see $\Theta^{[2]}$ is
  non-singular everywhere, Especially $\kappa^{-1}(\Theta^{[2]})$ is
  non-singular at the origin $(0,0,0)$.  On the other hand, let $u$ and
  $v$ be two variables on $\kappa^{-1}(\Theta^{[1]})$.  Then we have an
  expansion with respect to $v_3$:
 \[
 0=\sigma(u+v)=\sigma_3(u)v_3+\tfrac12\,\left(\sigma_2(u) +
   \sigma_{33}(u)\right){v_3}^2+\cdots,
 \]
 where $\sigma_i= \partial \sigma/\partial u_i$, etc.  Hence
 \[
 \sigma_3(u)=0  \qquad \sigma_2(u)+\sigma_{33}(u)=0.
 \]
Again by expansion
 \[
 0=\sigma_3(u)=\sigma_{33}(0)v_3+\cdots,
 \]
we see that
 \[
 \sigma_{33}(0)=0.
 \]
In summary, 
 \[
 \sigma_3(0)=\sigma_2(0)=0,
 \]
so from the above arguments and (\ref{L2.14.3}), we must have
 \[
 \sigma_1(0)\neq 0. 
 \]
Hence the $\sigma$-expansion must be odd.
\vskip 5pt
\noindent{\bf Step 2.}  Next we show that the terms of weight less
  than 5 vanish and $C_5(u)$ is non-trivial.  
 
  We write all the possible odd terms up to and including terms of
  weight $5$.  Using the first two equations in \ref{L4index}, we can
  show that the coefficients of the terms of weight four and less are
  zero, and that the coefficients of weight $5$ are given by those in
  $C_5$ up to multiplication by a constant.  We know from Step 1 that
  this constant is non-zero and we insert this constant into the
  $\varepsilon$.
%%%%%%%%%%%%%%%%%%%%%%%%%%%%%%%%%%%
\vskip 5pt
\noindent{\bf Step 3.}  We now calculate the coefficents $C_i$, $i>5$.
The proof of this step is by construction (with heavy use of Maple)
using the PDEs given in Lemma \ref{L4index}.  We expand
$\sigma(u_1,u_2,u_3)$ in a Taylor series with undetermined
coefficients, keeping only odd terms.  We do {\em not} assume that the
coefficients of the expansion are polynomial in the $\mu_i$, only that
they are independent of the $u_i$.  We then insert the expansion into
the 4-index PDEs for the $\wp$, and truncate to successive orders in
the weights of the $u_i$.  These give a series of linear equations for
the coefficients, and be using a sufficient number of the PDEs we can
always find unique solutions, as listed above.  We have carried out
this calculation down to $C_{18}$.  We have omitted the details of the
expressions for $C_{13}, \dots,C_{18}$, as they are rather lengthy,
but these are available at \cite{Weier34}.
%%%%%%%%%%%%%%%%%%%%%%%%%%%%%%%%%%%
\vskip 5pt
\noindent{\bf Step 4.}  Now consider the general term in the expansion.  Set
  \[
  A\,{u_1}^p{u_2}^q{u_3}^r, \quad
  A\in\mathbf{Q}(\mu_i)
  \]
  to be the lowest weight unknown term.  Since we have already shown by
  construction that the coefficients for all weights down to $-29$
  with respect to $u_j$s are polynomials, we may assume that
  $p+q+r\geqq 4$.  We consider the set $ (\sharp)$ of quadratic
  equations in $\sigma(u)$ and its (higher) derivatives obtained from
  the above, by multiplying the equations in \ref{L4index} by
  $\sigma(u)^2$.  We take an equation
  \begin{equation}
     \sigma(u)^2\,Q_{ijk\ell}(u)=\cdots\lbl{choosen} 
  \end{equation}
  from $(\sharp)$ such that $ {u_1}^p{u_2}^q{u_3}^r$ is divisible by
  $u_i u_j u_ku_{\ell}$.  We have at least one such equation.
  Differentiating (\ref{choosen}), we have an equation of the form
  \begin{equation}
    \sigma(u)\,
    \left(\frac{\partial^{p+q+r}\,\sigma}{\partial {u_1}^p\, 
        \partial {u_1}^q\,\partial {u_1}^r}\right)(u)
    +\cdots =0 \lbl{induction-key} 
  \end{equation}
such that all terms are polynomial of $\sigma(u)$ and its higher derivatives 
and such that
\[
\left(\frac{\partial^{p+q+r}\,\sigma}{\partial {u_1}^p\, 
\partial {u_1}^q\,\partial {u_1}^r}\right)(u)
\]
is the highest derivative in (\ref{induction-key}).  By looking at the
coefficient of the term $u_1$, we have a linear equation of the form
  \[
  A+\cdots=0 
  \]
over $\mathbf{Q}[\mu_1,\cdots,\mu_{12}]$.  Since the other terms
except $A$ in the above equation come from terms of $\sigma(u)$
whose weight is less than weight of ${u_1}^p{u_2}^q{u_3}^r$, 
we see $A$ is a polynomial in the $\mu_j$s by the induction hypothesis.
\end{proof}

\begin{remark}
(1) In Theorem \ref{L6.1}, the constant $\varepsilon$ might be unity, 
another 8th root of 1, or some other constant.  We have not been able to 
narrow down this result.  If the case $\varepsilon=1$ is true, then the 
determination of $\varepsilon$ reduces to the choice of roots in 
(\ref{discriminant}) and  (\ref{sigma-const}).  
The remaining results in this paper do not depend on this 
choice, or on the possibility that $\varepsilon \ne
1$.\lbl{varepsilon} 
\newline
(2) The weight of $\sigma(u)$ is inferred from  (\ref{sigma-const}) 
since the weight of $|\omega'|$ is $5+2+1$ and the conjectured weight 
of $D$ is $72$.  
The weight of the terms in the exponentials are all $0$ 
and the weight of $c$ is $72/8-(5+2+1)/2=5$ and coincides with 
the terms in the expansion of \ref{L6.1} if the weight 
of $\varepsilon$ is $0$.
\end{remark}

We shall need later on the following special property of the
$\sigma$-function in the purely trigonal case:
\begin{lemma} % 5.3
  The $\sigma$ function associated with the purely trigonal curve {\rm
    (\ref{eq5.1})} satisfies $ \sigma([-\zeta]u)=-\zeta\sigma(u)$ for
  $u \in \mathbb{C}^3$ under the notation {\rm(\ref{eq5.1})}.
  \lbl{sigsym}
\end{lemma}

\begin{proof} 
Since  $\Lambda$  is stable under the action of  $[\zeta]$  and  $[-1]$, 
we can check the statement by Lemma~\ref{L2.14} and Remark~\ref{R2.15}. 
\end{proof}

%%%%%%%%%%%%%%%%%%%%%%%%%%%%%%%%%%%%%%%%%%%%%%%%%%%%%%%%%%%%%%%%%%%%%%%%%%%%
\section{Basis of the space  $\varGamma(J, \mathcal{O}(n\Theta^{[2]}))$}
%%% Section 5 %%%
%%%%%%%%%%%%%%%%%%%%%%%%%%%%%%%%%%%%%%%%%%%%%%%%%%%%%%%%%%%%%%%%%%%%%%%%%%%%%
\label{BasisG}
For notational simplicity, we denote
\begin{equation}
   \partial_j=\tfrac{\partial}{\partial u_j}.
   \lbl{eq4.1}
\end{equation}
We also define
\begin{equation}
   \wp^{[ij]}=\text{the determinant of the $(i,j)$-(complementary) minor of 
$[\wp_{ij}]_{3\times 3}$}. \lbl{eq3.2b}
\end{equation}
We have explicit bases of the vector spaces
$\varGamma(J,\mathcal{O} (2\Theta^{[2]}))$ and
$\varGamma(J,\mathcal{O}(3\Theta^{[2]}))$ as follows (see also
\cite{cn06}, Example in Section 9):

\begin{lemma}% 4.2
We have the following\,{\rm :}\lbl{L4.2}
\begin{equation*}
\begin{aligned}
\varGamma(J, \mathcal{O}(2\Theta^{[2]}))&
  =     \mathbb{C} 1
   \oplus\mathbb{C} \wp_{11}
   \oplus\mathbb{C} \wp_{12}
   \oplus\mathbb{C} \wp_{13}
   \oplus\mathbb{C} \wp_{22}
   \oplus\mathbb{C} \wp_{23}
   \oplus\mathbb{C} \wp_{33}
   \oplus\mathbb{C} Q_{1333}, %\lbl{L4.2.1}
\\
\varGamma(J,\mathcal{O}(3\Theta^{[2]}))&= 
      \varGamma(J,\mathcal{O}(2\Theta^{[2]}))
   \oplus\mathbb{C} \wp_{111}
   \oplus\mathbb{C} \wp_{112}
   \oplus\mathbb{C} \wp_{113}
   \oplus\mathbb{C} \wp_{122}
   \oplus\mathbb{C} \wp_{123} \\
   & \quad
   \oplus\mathbb{C} \wp_{133} 
   \oplus\mathbb{C} \wp_{222}
   \oplus\mathbb{C} \wp_{223}
   \oplus\mathbb{C} \wp_{233}
   \oplus\mathbb{C} \wp_{333}
   \oplus\mathbb{C} \partial_1Q_{1333}
   \oplus\mathbb{C} \partial_2Q_{1333}\\
   & \quad
   \oplus\mathbb{C} \partial_3Q_{1333} 
   \oplus\mathbb{C} \wp^{[11]}
   \oplus\mathbb{C} \wp^{[12]}
   \oplus\mathbb{C} \wp^{[13]}
   \oplus\mathbb{C} \wp^{[22]}
   \oplus\mathbb{C} \wp^{[23]}
   \oplus\mathbb{C} \wp^{[33]}.
   %\lbl{L4.2.2}
  \end{aligned}
 \end{equation*}
\end{lemma}

\begin{proof} 
  We know the dimensions of the spaces above are $2^3=8$ and $3^3=27$,
  respectively by the Riemann-Roch theorem for Abelian varieties (see
  for example, \cite{mu85}, (pp.150--155), \cite{la82}, (p.99, Th.\
  4.1).  Moreover, (\ref{L2.14.3}) shows that the functions in the
  right hand sides belong to the spaces of the left hand sides,
  respectively.  For the space
  $\varGamma(J,\mathcal{O}(2\Theta^{[2]}))$, $\wp_{ij}$ and
  $Q_{ijk\ell}$ become the basis of the space from Definition
  \ref{D4.1}, Lemma \ref{L2.14}, and the arguments in the previous
  section.  However these are not all linearly independent, since
  there are connecting relations, such as those given in Lemma
  \ref{L4index},
  %\footnote{It should be noted that (\ref{kl1}) does not generate any such relations.}
  and the number of these relations
  is greater than the dimension of the space. Thus the problem is
  reduced to picking the linearly independent bases as a function
  space.  It is obvious that such independence does not depend upon
  the coefficients of curve by considering these expansions around the
  origin of $\mathbb C^3$.  Hence by multiplying by $\sigma(u)^2$ from
  the right hand side with respect to $u_1$, $u_2$, $u_3$, and after
  putting all the $\mu_j$ equal to zero, we see the functions of the
  right hand side are linearly independent.  The authors used a
  computer to check this.  Similarly, for the space
  $\varGamma(J,\mathcal{O}(3\Theta^{[2]}))$, the $27$ functions
  obtained by multiplying by $\sigma(u)^3$ from the right hand side
  are checked to be linearly independent by using a computer,
  expanding the given functions in the Abelian variables (cf.\ Theorem
  \ref{L6.1}) to a sufficiently high power that independence is
  checked.  We also see both decompositions in Lemma \ref{L4.2} in
  Example in Section 9 of \cite{cn06}.
\end{proof}

%%%%%%%%%%%%%%%%%%%%%%%%%%%%%%%%%%%%%%%%%%%%%%%%%%%%%%%%%%%%%%%%%%%%%%%%%
\section{The first main addition theorem} %%%% Section 9 %%%%%%%%%%%%
%%%%%%%%%%%%%%%%%%%%%%%%%%%%%%%%%%%%%%%%%%%%%%%%%%%%%%%%%%%%
\label{Add1}
\begin{theorem} % 7.1
  The $\sigma$-function associated with {\rm (\ref{eq5.1})} satisfies
  the following addition formula on $J\times J$\,{\rm :}\lbl{T7.1}
\begin{equation}
  \begin{split}
    -\frac{\sigma(u+v)\sigma(u-v)}{\sigma(u)^2\sigma(v)^2}&=
    \wp_{11}(u) -\wp_{11}(v) +\wp_{12}(u)\wp_{23}(v)
    -\wp_{12}(v)\wp_{23}(u)\\
    & +\wp_{13}(u)\wp_{22}(v) -\wp_{13}(v) \wp_{22}(u)
    +\tfrac13\left(\wp_{33}(u)Q_{1333}(v) -\wp_{33}(v) Q_{1333}(u)
    \right)
    \\
    & -\tfrac13\mu_1 \left( \wp_{12}(u) \wp_{33}(v)
      -\wp_{12}(v) \wp_{33}(u) \right) -\mu_1 \left(
      \wp_{13}(u)\wp_{23}(v) -\wp_{13}(v) \wp_{23}(u)
    \right)\\
    & +\tfrac13 \left({\mu_1}^2 -\mu_2 \right) \left(
      \wp_{13}(u) \wp_{33}(v) -\wp_{13}(v)\wp_{33}(u) \right)
    +\tfrac13 \mu_8 \left(\wp_{33}(u)-\wp_{33}(v)\right)
  \end{split}
\end{equation}
\end{theorem}

\begin{proof}
  Firstly, we notice that the left hand side is an odd function with
  respect to $(u,v)\mapsto ([-1]u,[-1]v)$, and that it has poles of
  order 2 along $(\Theta^{[2]}\times J)\cup(J\times\Theta^{[2]})$ but
  nowhere else.  Moreover it is of weight $-10$.  
  %{\color{blue}
  Therefore, by Lemma \ref{L4.2}, the left hand side is expressed by 
  a finite sum of the form 
  \begin{equation}
  \sum_jA_j\,\big(X_j(u)Y_j(v)-X_j(v)Y_j(u)\big), \lbl{rhs}
  \end{equation}
  where the $A_j$ are rational functions of the $\mu_i$s with homogeneous
  weight, and the $X_j$ and $Y_j$ are functions chosen from the right hand
  side of the first equality in Lemma \ref{L4.2}.  We claim that all
  the $A_j$ are polynomial in the $\mu_i$s.  Suppose all the $A_j$s
  are reduced fractional expressions, and at least one of the $A_j$s
  is not a polynomial.  Take the least common multiple $B$ of all the
  denominators of the $A_j$s.  Note that there is a set of special
  values of the $\mu_i$s such that $B$ vanishes and the numerator of
  at least one $A_j$ does not vanish.  After multiplying the equation
  $\mbox{\lq\lq lhs"}{=}$\,(\ref{rhs}) by $B\,\sigma(u)^2\sigma(v)^2$,
  and taking the $\mu_i$s to be such a zero of $B$, we have a
  contradiction, by using the linear independency of Lemma \ref{L4.2} twice
  with respect to the variables $u$ and $v$ for the corresponding
  curve of (\ref{eq1.1}).  Hence, all the $A_j$ must be polynomials.
  %}
%  If we regard the right hand side as a rational function of $\sigma(u)$, $\sigma(v)$,
%  $\sigma_i(u)$, $\sigma_i(v)$, $\sigma_{ij}(u)$, $\sigma_{ij}(v)$,
%  $\cdots$, and consider the two sides multiplied by
%  $\sigma(u)^2\sigma(v)^2$, 
  % which are homogeneous power series of
  % weight $30$ ($\in\mathbb{Q}[\mu_i] [[u_1,u_2, u_3]]$), 
  Hence, we see that the desired right hand side must be expressed by using
  constants $a,b,c,d,e,f,g_1,g_2,h_1,h_2,i_1,i_2,j,k_1,k_2,k_3$ which are 
  polynomials in  $\mu_i$s and independent of the $u_i$ and $v_i$, as follows:
\begin{equation}
   \begin{aligned}
          & a\,[\wp_{11}(u) - \wp_{11}(v)]
           +b\,[\wp_{12}(u)\wp_{23}(v) - \wp_{12}(v)\wp_{23}(u)] 
           +c\,[\wp_{13}(u)\wp_{22}(v) - \wp_{13}(v)\wp_{22}(u)] \lbl{eq7.2}\\
     &\quad+d\,[Q_{1333}(u)\wp_{33}(v) - Q_{1333}(v)\wp_{33}(u)]
           +e\mu_1[\wp_{12}(u)\wp_{33}(v)-\wp_{12}(v)\wp_{33}(u)] \\
     &\quad+f[\wp_{13}(u)\wp_{23}(v)-\wp_{13}(v)\wp_{23}(u)]
           +g_1[\wp_{13}(u)\wp_{33}(v)-\wp_{13}(v)\wp_{33}(u)]\\
     & \quad
           +g_2[Q_{1333}(u) - Q_{1333}(v)]
      +h_1[\wp_{23}(u)\wp_{22}(v)-\wp_{23}(v)\wp_{22}(u)]
     +h_2[\wp_{12}(u)-\wp_{12}(v)]\\ & \quad
     +i_1[\wp_{22}(u)\wp_{33}(v)-\wp_{22}(v)\wp_{33}(u)]
     +i_2[\wp_{13}(u)-\wp_{13}(v)]+j[\wp_{23}(u)\wp_{33}(v)-\wp_{23}(v)\wp_{33}(u)]\\ 
     & \quad +k_1[\wp_{22}(u)-\wp_{22}(v)] 
             +k_2[\wp_{23}(u)-\wp_{23}(v)]  +k_3[\wp_{33}(u)-\wp_{33}(v)].
   \end{aligned}
\end{equation}
We find by computer using Maple, on substituting the expansion
(\ref{L6.1}) up to $C_{13}$ terms of $\sigma(u)$ into (\ref{eq7.2}),
and truncating up to weight 18 in the $u_i$ and $v_i$, that
\begin{equation}
 \begin{split}
   & a=b=c=-1, \quad d=\tfrac13, \quad
   e=-\tfrac13\mu_1,\quad f=-\mu_1,\quad
   g_1=\tfrac13 (\mu_1^2-\mu_2),\\
   & \quad g_2=h_1=h_2=i_1=i_2=j=k_1=k_2=0, \quad  k_3=\tfrac13 \mu_8. \lbl{eq7.3}
\end{split} 
\end{equation}
as asserted. In the Maple calculation, it is not necessary to assume
the polynomial nature of the coefficients as functions of the $\mu_j$.
\end{proof}

\begin{remark} % 7.4 
By applying \lbl{R7.4}
\begin{equation}
   \frac12\frac{\partial}{\partial u_i}
   \left(
   \frac{\partial}{\partial u_j}
  +\frac{\partial}{\partial v_j}
   \right)\log
\end{equation}
to  \ref{T7.1}, we have $-\wp_{ij}(u+v)+\wp_{ij}(u)$ 
from the left hand side, and have a rational expression
of several $\wp_{ij\cdots\ell}(u)$s and $\wp_{ij\cdots\ell}(v)$s on
the right hand side.  Hence, we have an algebraic addition formulae for
$\wp_{ij}(u)$s.
\end{remark}
\begin{remark}
By putting $v=u-(\delta,0,0)$ and letting $\delta \rightarrow 0$, we
can get a ``double-angle'' $\sigma$-formula \lbl{R9.3}
\begin{equation}
\begin{split}
  \frac{\sigma(2u)}{\sigma(u)^4}& = -\wp_{111}(u)-\wp_{112}\wp_{23}+
  \wp_{12}(u)\wp_{123}(u)-\wp_{113}(u)\wp_{22}(u)+\wp_{13}(u)\wp_{122}(u)\\
  &\quad -\tfrac13\wp_{133}(u) Q_{1333}(u) 
         +\tfrac13\wp_{33}(u)\tfrac{\partial}{\partial u_1}Q_{1333}(u)
   +\tfrac13\mu_1\big(\wp_{112}(u)\wp_{33}(u)-\wp_{12}(u)\wp_{133}(u)\big)\\
  &\quad +\mu_1\left(\wp_{113}(u)\wp_{23}-\wp_{13}(u)\wp_{123}\right)
         -\tfrac13\big({\mu_1}^2 -\mu_2\big)
                  \big( \wp_{113}(u) \wp_{33}(u)-\wp_{13}(u)\wp_{133}(u)\big)\\
  &\quad 
         -\tfrac13 \mu_8 \wp_{133}(u).
\end{split}
\end{equation}
In the case of the elliptic curve, the corresponding relation is
$\sigma(2u)=-\wp'(u)\sigma^4(u)$, whilst the corresponding formula for
the hyperelliptic genus two curve is given in \cite{ba07}, p. 129.
\end{remark}

%%%%%%%%%%%%%%%%%%%%%%%%%%%%%%%%%%%%%%%%%%%%%%%%%%%%%%%%%%%%%%%%%%%
\section{The second main addition theorem} %%%%%%% Section 10 %%%
%%%%%%%%%%%%%%%%%%%%%%%%%%%%%%%%%%%%%%%%%%%%%%%%%%%%%%%%%%%%%%%%%%%
\label{Add2}
The second main addition result applies only in the purely trigonal
case (\ref{eq5.1}), using the results of Lemma \ref{sigsym}.  The
formula is as follows:
\begin{theorem} % 8.1
The $\sigma$-function associated with {\rm (\ref{eq5.1})} satisfies the
following addition formula on $J\times J $\,{\rm :}\lbl{T8.1} 
\begin{equation}
    \frac{\sigma(u+v)\sigma(u+[\zeta]v)\sigma(u+[\zeta^2]v)}
{\sigma(u)^3\sigma(v)^3}
     = R(u,v)+R(v,u)\lbl{zeta-add},
\end{equation}
where
\begin{equation*}\begin{split}
R(u,v) 
    &= \quad-\tfrac13\wp_{13}(u)\partial_3Q_{1333}(v) 
            -\tfrac34\wp_{23}(u)\wp_{112}(v)  -\tfrac12 \wp_{111}(u) 
            +\tfrac14\wp_{122}(u)\wp^{[11]}(v)  \\
    & \quad -\tfrac14\wp_{222}(u)\wp^{[12]}(v) 
            +\tfrac1{12}\partial_3Q_{1333}(u)\wp^{[11]}(v) 
            +\tfrac12\wp_{333}(u)\wp^{[22]}(v) 
           -\tfrac14 \mu_3 \wp_{333}(u)\wp^{[12]}(v)   \\
    & \quad+\tfrac12 \mu_6 \wp_{13}(u)\wp_{333}(v) 
        -\tfrac14\mu_9\wp_{23}(u)\wp_{333}(v)
    -\tfrac12  \mu_{12} \wp_{333}(u). 
\end{split}\end{equation*}
\end{theorem}

\begin{proof} % of  8.1. 
Our goal is to express 
\begin{equation}
   \frac{\sigma(u+v)\sigma(u+[\zeta]v)\sigma(u+[\zeta^2]v)}
        {\sigma(u)^3\sigma(v)^3}
   \lbl{eq8.3}
\end{equation}
using several $\wp$ functions.
  %{\color{blue}
Because (\ref{eq8.3}) belongs to $\varGamma(J\times J,
\mathcal{O}(3((\Theta^{[2]}\times J)\cup(J\times \Theta^{[2]}))))$, a
similar argument to that at the beginning of the proof of Th.\
\ref{T7.1} shows that it must be a finite sum of multi-linear forms of
the 27 functions in Lemma \ref{L4.2}, namely, of the form
\begin{equation}
   \sum_{j}^{\text{finite sum}}\hskip -3pt C_j\,X_j(u)Y_j(v),
   \lbl{eq8.4}
\end{equation}
where $X_j$ and $Y_j$ are any of the functions appearing in the right
hand side of the description of $\varGamma(J,\mathcal{O}(3\Theta^{[2]})))$
in Lemma \ref{L4.2}, and the $C_j$ are polynomial in $\mu_i$s.  
%}
Moreover, (\ref{eq8.3}) has the following properties:
\begin{enumerate}
\item[L1.] As a function on  $J\times J$, its weight is
   $(-5)\times 3=-15$;
\item[L2.] It is invariant under  $u\mapsto [\zeta]u$
   (resp. $v\mapsto [\zeta]v$);
\item[L3.] It has a pole of order $3$ on  $(\Theta^{[2]}\times
   J)\cup(J\times \Theta^{[2]})$;
\item[L4.] It is invariant under the exchange  $u\leftrightarrow v$
  (by Lemma \ref{sigsym}). 
\end{enumerate}
Hence, (\ref{eq8.4}) has the same properties.  Thus, we may consider
only the functions in our basis of $\varGamma(J, \mathcal{O}
(3\Theta^{[2]}))$ that have the following corresponding properties:
\begin{itemize}
\item[R1.] The weight is greater than or equal to  $(-5)\times
   3=-15$;
\item[R2.] They are invariant under  $u\mapsto [\zeta]u$;
\item[R3.] They have poles of order at most $3$ on  $\Theta^{[2]}$.
\end{itemize}
There are 12 such functions and they are listed as follows:
\begin{align*}
  & 1, && \wp_{13} \quad (\text{weight}=-6),
  && \wp_{23} \quad (\text{weight}=-3),   \\
  & \wp_{111} \quad (\text{weight}=-15), && \wp_{112} \quad
  (\text{weight}=-12),
  && \wp_{122} \quad (\text{weight}=-9),  \\
  & \wp_{222} \quad (\text{weight}=-6), && \wp_{333} \quad
  (\text{weight}=-3),
  && \wp^{[22]} \quad (\text{weight}=-12), \\
  & \wp^{[12]} \quad (\text{weight}=-9), && \wp^{[11]} \quad
  (\text{weight}=-6), && \end{align*}
    \[ \partial_3 Q_{1333}=
         -6(\wp_{13}\wp_{333}-\wp_{133}\wp_{33})-3\wp_{122}
   \quad (\text{weight}=-9),
\]
and the $\wp^{[ij]}$ are defined in (\ref{eq3.2b}).
Here the last equality is given by
cross-differentiation from $\partial_1 Q_{3333}$ using the first of
the relations in Lemma \ref{L4index} with $\mu_1=\mu_2=\mu_4=0$.
Since (\ref{eq8.3}) is an even function, it must be of the form
%\newpage
\begin{equation}
\frac{\sigma(u+v)\sigma(u+[\zeta]v)\sigma(u+[\zeta^2]v)}
{\sigma(u)^3\sigma(v)^3}
     = \tilde R(u,v)+ \tilde R(v,u),\lbl{eq8.6}
\end{equation}
where  
\begin{align*}
  \tilde R(u,v) &= a_1\wp_{13}(u)\wp_{122}(v)
           +a_2\wp_{13}(u)\,\partial_3Q_{1333}(v)
           +a_3\wp_{23}(u)\wp_{112}(v)
           +a_4\wp_{111}(u) \\
  & \quad  +a_5\wp_{122}(u)\wp^{[11]}(v)
           +a_6\wp_{222}(u)\wp^{[12]}(v)
           +a_7\,\partial_3Q_{1333}(u)\wp^{[11]}(v)
           +a_8\wp_{333}(u)\wp^{[22]}(v) \\
  & \quad  +b_1\wp_{13}(u)\wp_{222}(v)
           +b_2\wp_{23}(u)\wp_{122}(v)
           +b_3\wp_{23}(u)\,\partial_3Q_{1333}(v)   +b_4\wp_{112}(u)\\
  & \quad
           +b_5\wp_{222}(u)\wp^{[11]}(v) 
           +b_6\wp_{333}(u)\wp^{[12]}(v)  +c_1\wp_{13}(u)\wp_{333}(v)
          +c_2\wp_{23}(u)\wp_{222}(v)\\
  & \quad
          +c_3\wp_{122}(u)  +c_4\wp_{333}(u)\wp^{[11]}(v)
          +c_5\,\partial_3Q_{1333}(u) +d_1\wp_{23}(u)\wp_{333}(v)
          +d_2\wp_{222}(u) \\
  & \quad +e_1\wp_{333}(u). 
\end{align*}
By substituting (\ref{L6.1}) into (\ref{eq8.6}), and comparing
coefficients of different mononomials in $u_i,v_j$, we can find the
constants $a_1$, $\cdots$, $e_1$ depending on the $\mu_k$s.  Again, in
this lengthy Maple calculation, it is not necessary to assume the
coefficients are polynomial in the $\mu_i$.
% To avoid excessive
%   computing time, we divide up this calculation into separate terms
%   involving independent powers of the $\mu_j$.
\end{proof}

\begin{remark} % 8.2
By applying 
\begin{equation}
   \frac13\left(
   \frac{\partial^2}{\partial u_i\partial u_j}
  +\frac{\partial^2}{\partial u_i\partial v_j}
  +\frac{\partial^2}{\partial v_i\partial v_j}
   \right)\log
\end{equation}
   to (\ref{zeta-add}), we obtain  algebraic addition
   formulae for standard Abelian functions,
which would be interesting to compare with those of Remark (\ref{R7.4}).
\end{remark}

\begin{remark} By putting $v=-u+(\delta,0,0)$ into (\ref{zeta-add}),
  dividing through by $\delta$ and letting $\delta \rightarrow 0$, we
can get an unusual ``shifted'' $\sigma$-formula of the form \lbl{R10.2}
\begin{align}
  &-\frac{\sigma(u-[\zeta]u)\sigma(u-[\zeta^2]u)}{\sigma(u)^6}
  =\sum_{i=1}^{12} c_i \left[g_i(u) \partial_1f_i(u)-f_i(u)\partial_1
    g_i(u)\right],\lbl{shifted}
\end{align}
where the $f_i$ and the $g_i$ are the even and odd derivative
components respectively of the formula in (\ref{zeta-add}), i.e.\ as
given in the following table
\begin{equation*}
\setlength{\extrarowheight}{2pt}
\begin{array}{c|c|c||c|c|c}c_i & f_i & g_i& c_i & f_i & g_i \\\hline
 \tfrac12          & \wp_{13}(u) & \wp_{122}(u)          &  -\tfrac13      &  \wp_{13}(u) & \partial_3Q_{1333}(u) \\
-\tfrac34          & \wp_{23}(u) & \wp_{112}(u)          &  -\tfrac12      &        1     & \wp_{111}(u) \\
 \tfrac14          &\wp^{[11]}(u)& \wp_{122}(u)          &  -\tfrac14      & \wp^{[12]}(u)& \wp_{222}(u)  \\
 \tfrac1{12}       &\wp^{[11]}(v)& \partial_3Q_{1333}(u) &   \tfrac12      & \wp^{[22]}(u)& \wp_{333}(u)\\
-\tfrac14\mu_3 &\wp^{[12]}(u)& \wp_{333}(u)          &\tfrac12\mu_6& \wp_{13}(u)  & \wp_{333}(u) \\
-\tfrac14\mu_9 & \wp_{23}(u) & \wp_{333}(v)          &-\tfrac12  \mu_{12} 
                                                                           &        1     & \wp_{333}(u)
\setlength{\extrarowheight}{0pt}
\end{array}\end{equation*}
\end{remark}
\begin{remark}
In the general elliptic case, there appears to be no formulae
corresponding to (\ref{zeta-add}) and (\ref{shifted}).  However for the
specialized {\em equianharmonic case}, where $\wp$ satisfies
\[
(\wp')^2 = 4 \wp^3 - g_3, 
\]
it is straightforward to show that
\[
\frac{\sigma(u+v)\sigma(u+\zeta v)\sigma(u+\zeta^2 v)}
{\sigma^3(u)\sigma^3(v)} 
   = -\tfrac12 (\wp'(u) + \wp'(v)),
\]
and
\[
\frac{\sigma\left((1-\zeta)u\right)\sigma\left((1-\zeta^2)u\right)}
{\sigma^6(u)} = 3 \wp^2(u).
\]
These seem to be just the first of a family of multi-term addition
formulae on special curves
with automorphisms, which will be discussed in
more detail elsewhere.
\end{remark}
\section*{Acknowledgements}
\addcontentsline{toc}{section}{Acknowledgements}
This paper was started during a visit by the authors to Tokyo
Metropolitan University in 2005, supported by JSPS grant 16540002.  We
would like to express our thanks to Prof.\ M.\ Guest of TMU who helped
organize this visit.  The work continued during a visit by VZE to
Heriot-Watt University under the support of the Royal Society.
Further work was done whilst three of the authors (JCE, VZE, and EP)
were attending the programme in Nonlinear Waves at the Mittag-Leffler
Institute in Stockholm in 2005, and we would like to thank Professor
H.\ Holden of Trondheim and the Royal Swedish Academy of Sciences for
making this possible [EP, being then on leave from Boston University,
is grateful for NSA grant MDA904-03-1-0119 which supported her
doctoral students who were performing related research].  The authors
are also grateful for a number of useful discussions with
Prof.~A.~Nakayashiki, and Drs.~John Gibbons and Sadie Baldwin.  In
particular we are grateful to John Gibbons for pointing out the
possibility of the relations described in Remarks \ref{R9.3} and
\ref{R10.2}.  We are grateful to Mr Matthew England for pointing out a
number of typos in various versions of this manuscript.  Some of the
calculations described in this paper were carried out using
Distributed Maple \cite{smb03}, and we are grateful to the author of
this package, Professor Wolfgang Schreiner of RISC-Linz, for help and
advice. 
%{\color{blue} 
Finally, we would like to express special thanks to the referees for
constructive suggestions to improve the paper, in particular for
pointing out some crucial gaps in the main theorems and for giving
hints how to fill them. %}

%\newpage

\setcounter{section}{0}
\setcounter{equation}{0}
\renewcommand{\thesection}{\Alph{section}}
\section*{Appendix A:  The fundamental bi-differential}
\addtocounter{section}{1}
\addcontentsline{toc}{section}{Appendix A:  The fundamental bi-differential}

We write the polynomial $f(x,y)$ in (\ref{eq1.1}) that defines the
trigonal curve $C$ as
\begin{equation}
f(x,y)=y^3+p(x)y^2+q(x)y-r(x)\lbl{trigonal}
\end{equation} 
with 
\begin{equation*}
p(x)=\mu_1x+\mu_4, \quad \
q(x)=\mu_2x^2+\mu_5x+\mu_8, \quad \
r(x)=x^4 +\mu_3 x^3 +\mu_6 x^2 +\mu_9 x +\mu_{12}.
\end{equation*}
We describe explicitly the fundamental non-normalized bi-differential
(Klein's fundamental 2-form of the second kind) $\Sigma((x,y),(z,w))$
in (\ref{eq2.3}) of the curve for $(x,y)$, $(z,w)$ in $C$ defined by
$f(x,y)=0$.

Following the scheme described in \cite{ba97} and applied to trigonal
curves in \cite{eel00}, \cite{bel00} and the present paper, one can
realize $\Sigma((x,y),(z,w))$ explicitly as
\begin{equation} 
  \Omega((x,y),(z,w))=\frac{F((x,y),(z,w))dxdz}{(x-z)^2 f_y(x,y)f_w(z,w)}   
\end{equation}
with the polynomial $F((x,y);(z,w))$ given by the formula 
\begin{align}\begin{split}
F\big(&(x,y),(z,w)\big)=(wy+Q(x,z))(wy+Q(z,x))\lbl{polar}\\
&+w\left(w \left[\frac{f(x,y)}{y}\right]_y +T(x,z)\right)
 +y\left(y \left[\frac{f(z,w)}{w}\right]_w +T(z,x)\right)
 -F_0(x,z)\end{split}
\end{align}  
with 
\begin{align}
 \begin{split}
    Q(x,z)&=( \mu_1^{2}-\mu_2) xz
   + (2\,\mu_1\mu_4-\mu_5) x-\mu_8+\mu_4^{2}\lbl{QT}\\
    T(x,z)&=3\mu_{12}+ ( z+2x ) \mu_9+x ( x+2\,z )
    \mu_6+3\mu_3 x^2 z+ p(z) q(x) + x^2 z^{2}+2\,x^3z.
\end{split}
\end{align} 
The term $F_0(x,y)$ vanishes at $\mu_1=\mu_4=0$ and is given by 
\begin{align*}
F_0(x,y)&=c_{32}(x+z)x^2z^2+c_{22}x^2z^2+c_{21}(x+z)xz+c_{11}xz+c_{10}(x+z)+c_{00}, \\
c_{32} &=-\mu_1, \qquad c_{22}=-2\mu_4-2{\mu_1}^2\mu_2+{\mu_1}^4+2\mu_3\mu_1, \\
c_{21} & =\mu_6\mu_1-2\,\mu_1\mu_4\mu_2+\mu_3\mu_4-\mu_5{\mu_1}^2+2{\mu_1}^3\mu_4,\\
c_{11} & = 2\left( 3{\mu_1}^2{\mu_4}^2+\mu_6\mu_4+\mu_9\mu_1-2\,\mu_1\mu_4\mu_5
          -{\mu_1}^2\mu_8-\mu_2{\mu_4}^2 \right),\\
c_{10} &=  -\mu_5{\mu_4}^2+\mu_1\mu_{12}+2\,\mu_1{\mu_4}^3-2\,\mu_1\mu_4\mu_8+\mu_3\mu_4,
\qquad c_{00} = \mu_4 \left( {\mu_4}^3+2\,\mu_{12}-2\,\mu_4\mu_8 \right). 
\end{align*}
We also remark that the expression (\ref{polar}) generalizes 
the Kleinian 2-polar previously derived in the hyperelliptic case
\cite{ba97}.

\section*{Appendix B: Quadratic three-index relations}
\addcontentsline{toc}{section}{Appendix B: Quadratic three-index relations}
A complete list of the known relations quadratic in three-index
$\wp_{ijk}$, up to weight $-15$, for the ``purely trigonal'' case is
given below.  Note that with care we can obtain an expression such
that the highest power on the r.h.s.\ is no more than cubic.  The
number in square brackets [\quad ] is the weight.  A fuller list for the
{\em general} (3,4) case is given at \cite{Weier34}. 
\begin{align*}
  {\wp_{333}}^2 & = {\wp_{23}}^2+4 \wp_{13} - 4\wp_{33} \wp_{22} 
   + 4{\wp_{33}}^3, \quad [-6]
  \\
  \wp_{233} \wp_{333} & = -\wp_{22}\wp_{23} + 4 \wp_{23}{\wp_{33}}^2 
   + 2\mu_3{\wp_{33}}^2 - 2 \wp_{12}, \quad [-7]
  \\
  {\wp_{233}}^2 & = 4 \wp_{33}{\wp_{23}}^2 + 4\mu_3\wp_{33}\wp_{23}
  + {\wp_{22}}^2 -\tfrac43 \wp_{1333} + 4\mu_6\wp_{33}  \nonumber
  +8\wp_{33}\wp_{13}, \quad [-8]
  \\
  \wp_{333}\wp_{223} & = 2 \wp_{33}{\wp_{23}}^2 + \mu_3\wp_{33} 
  \wp_{23} - 2{\wp_{22}}^2 + \tfrac23 \wp_{1333} + 2{\wp_{33}}^2 \wp_{22},
  \quad [-8]
  \\
  \wp_{223} \wp_{233} & = 2{\wp_{23}}^3 + 2\wp_{33}\wp_{23}\wp_{22}
   + 2\mu_9 + 4 \wp_{23}\wp_{13} + 2\mu_6\wp_{23}  
   + 2\mu_3\wp_{13} +2 \mu_3{\wp_{23}}^2 \\
& \qquad+ \mu_3\wp_{22}\wp_{33}, \quad
  [-9]  \\
  \wp_{222} \wp_{333} & = -\mu_3{\wp_{23}}^2 - 4\mu_3\wp_{13}
  + 4\mu_3\wp_{22}\wp_{33} - 2 {\wp_{23}}^3 - 8 \wp_{13}\wp_{23}
  + 6\wp_{33}\wp_{23}\wp_{22}  - 4\wp_{33}\wp_{12}, \ [-9] \\
  {\wp_{223}}^2 & = 4{\wp_{23}}^2\wp_{22} + 4 \wp_{11} + 4\wp_{22}\wp_{13}
  + 4\mu_6\wp_{22} - 4\wp_{23}\wp_{12} - 4\mu_3\wp_{12}\nonumber\\
  &\quad + {\mu_3}^2{\wp_{33}}^2 - 4\mu_6{\wp_{33}}^2 
  + \tfrac43\wp_{33}\wp_{1333} - 8 \wp_{13}{\wp_{33}}^2
  + 4\mu_3\wp_{23}\wp_{22}, \ [-10]  \\
  \wp_{133} \wp_{333} & = -2 \wp_{22} \wp_{13}+\wp_{12} \wp_{23}
  +\tfrac23 \wp_{33} \wp_{1333}, \quad [-10]  \\
  \wp_{233} \wp_{222} & = 2\mu_3\wp_{12} - \tfrac83\wp_{33}\wp_{1333}
  + 2\wp_{33}{\wp_{22}}^2 + 8\mu_6{\wp_{33}}^2 + 16{\wp_{33}}^2 \wp_{13}
  \nonumber\\
  & \quad -2{\mu_3}^2{\wp_{33}}^2 + 4\wp_{23}\wp_{12} 
  + \mu_3\wp_{23}\wp_{22} + 2{\wp_{23}}^2 \wp_{22}, \quad [-10]  \\
  \wp_{123} \wp_{333} & = 4\wp_{33}\wp_{23}\wp_{13} 
  + 2\mu_3 \wp_{33}\wp_{13} - 2\wp_{22}\wp_{12} - \tfrac13\wp_{23}\wp_{1333}
  +2{\wp_{33}}^2\wp_{12},\nonumber \quad [-11]   \end{align*}\begin{align*}
  \wp_{223} \wp_{222} & = 8 \wp_{33}\wp_{13}\wp_{23} 
  - \tfrac23\mu_3\wp_{1333} + 4\mu_9 \wp_{33} + 4\mu_6\wp_{33} \wp_{23} 
  + 4\mu_3\wp_{33}\wp_{13}\nonumber\\ 
  &\quad +2 \mu_3{\wp_{22}}^2 
  - \tfrac43\wp_{23}\wp_{1333} - {\mu_3}^2\wp_{33}\wp_{23}
  + 4 \wp_{23}{\wp_{22}}^2, \quad [-11]   \\
  \wp_{233} \wp_{133} & = 2\mu_9\wp_{33} + 2\mu_3\wp_{33}\wp_{13}
  + \tfrac23 \wp_{23}\wp_{1333} + \wp_{22}\wp_{12}, \quad [-11]   \\
  \wp_{123}\wp_{233} & = 2\wp_{33}\wp_{23}\wp_{12} 
  + 2\mu_3\wp_{33}\wp_{12} - 2\wp_{33}\wp_{11} 
  - \tfrac13\wp_{22}\wp_{1333} + 2\wp_{33}\wp_{22}\wp_{13} \nonumber \\
  &\quad + 2\mu_3\wp_{23}\wp_{13} + 2{\wp_{13}}^2 +2\mu_6\wp_{13}
  +2{\wp_{23}}^2\wp_{13} + 2 \mu_{12},\nonumber \quad [-12] \\
  \wp_{333}\wp_{122} & = - 2\wp_{33}\wp_{22}\wp_{13} 
  - \mu_3\wp_{23}\wp_{13} - 2\wp_{13}{\wp_{23}}^2 
  - 6{\wp_{13}}^2 - 2\mu_6\wp_{13} + \tfrac23\wp_{22}\wp_{1333}  \nonumber \\ 
  &\quad + 4\wp_{33}\wp_{23}\wp_{12} + 2\mu_3\wp_{33}\wp_{12}
  - 2 \wp_{33}\wp_{11} + 2 \mu_{12}+\mu_9 \wp_{23}, \quad  [-12]\nonumber \\ 
  \wp_{223} \wp_{133} & = 2\mu_3\wp_{23}\wp_{13} + 2{\wp_{23}}^2\wp_{13}
  + 2{\wp_{13}}^2 + 2\mu_6 \wp_{13} - \mu_3\wp_{33}\wp_{12} \nonumber \\ 
  & \quad - 2\wp_{33}\wp_{22}\wp_{13} + \tfrac23\wp_{22}\wp_{1333} 
  + 2\wp_{33}\wp_{11} +2\mu_{12}, \quad [-12] \\
%    \end{gaalign*}
%\begin{align*}
  {\wp_{222}}^2 & =  - 4\mu_3\wp_{33}\wp_{12} + 8\wp_{33}\wp_{11}
  - 4\wp_{22}\wp_{1333} + 24\wp_{33}\wp_{22}\wp_{13}\nonumber \\ &
  \quad + 4\mu_3\wp_{23}\wp_{13} - 8{\wp_{13}}^2 - 4 \mu_9\wp_{23}
  - 8\wp_{13}\mu_6 + 4{\wp_{22}}^3 + 4{\mu_3}^2\wp_{13} \nonumber \\ 
  & \quad + {\mu_3}^2{\wp_{23}}^2 - 4\mu_6{\wp_{23}}^2 
  - 4{\mu_3}^2 \wp_{33}\wp_{22} + 16\mu_6\wp_{33}\wp_{22} 
  - 8\mu_{12}, \quad [-12]  \\
  \wp_{223} \wp_{123} & = -2\wp_{23}\wp_{11} + 2\mu_9\wp_{22} 
  + 2\wp_{13}\wp_{23}\wp_{22} + 2\mu_3\wp_{23}\wp_{12}  
  + 2\wp_{12}{\wp_{23}}^2 \nonumber \\ 
  & \quad + \tfrac13\mu_3\wp_{33}\wp_{1333}
  - 2\mu_3\wp_{13}{\wp_{33}}^2 - 2\mu_9{\wp_{33}}^2 
  + 2\mu_3\wp_{22}\wp_{13}, \quad [-13]  \\
  \wp_{133}\wp_{222} & = 4\wp_{23}\wp_{22}\wp_{13}
  - \mu_3\wp_{23}\wp_{12} - 2{\wp_{23}}^2\wp_{12}
  - \tfrac23\mu_3\wp_{33}\wp_{1333} 
  + 4\mu_3\wp_{13}{\wp_{33}}^2  \nonumber \\ &
  \quad  + 4\mu_9{\wp_{33}}^2 + 2\mu_3\wp_{22}\wp_{13} 
  + 2\wp_{33}\wp_{22}\wp_{12}, \quad [-13]
  \\
  \wp_{122}\wp_{233} & = -\mu_9\wp_{22} + 4\wp_{13}\wp_{12}
  + 2\mu_6\wp_{12} - \tfrac23\mu_3\wp_{33}\wp_{1333}  
  + 2\wp_{33}\wp_{22}\wp_{12}  \nonumber \\ & 
  \quad + 4\mu_3{\wp_{33}}^2\wp_{13} + 4\mu_9{\wp_{33}}^2 
  + 2\mu_3\wp_{23}\wp_{12} + 2{\wp_{23}}^2\wp_{12} 
  - \mu_3\wp_{22}\wp_{13},  \nonumber \qquad [-13]  \\
  \wp_{333}\wp_{113} & = -2{\wp_{12}}^2 - \tfrac23\wp_{13}\wp_{1333}
  + 6\wp_{33}{\wp_{13}}^2 + 2\mu_6\wp_{33}\wp_{13} 
  + 2{\wp_{33}}^2\wp_{11} - 2\mu_{12}\wp_{33} \\
& \qquad
  - \mu_9\wp_{33}\wp_{23}, \quad [-14]  \\
  {\wp_{133}}^2 & = \tfrac43\wp_{13}\wp_{1333} + {\wp_{12}}^2 
  - 4\wp_{33}{\wp_{13}}^2 + 4\mu_{12}\wp_{33}, \quad [-14]  \\
  \wp_{223}\wp_{122} & = - 2\wp_{11}\wp_{22} + 4\wp_{23}\wp_{22}\wp_{12}
  + \tfrac43\wp_{13}\wp_{1333} + 2{\wp_{12}}^2 - 8\wp_{33}{\wp_{13}}^2 
  \\
& \qquad + \mu_3(2\wp_{22}\wp_{12} + 4\wp_{33}\wp_{23}\wp_{13}) -\tfrac23\mu_3\wp_{23}\wp_{1333} + {\mu_3}^2\wp_{33}\wp_{13}
  - \tfrac23\mu_6\wp_{1333}  \nonumber \\ & 
  \quad + 4\mu_9 \wp_{33}\wp_{23} 
  + (8\mu_{12} +\mu_3\mu_9)\wp_{33}, \quad [-14]  \\
\wp_{123}\wp_{222} & = 2{\wp_{22}}^2\wp_{13} 
  + 2\wp_{23}\wp_{22}\wp_{12} - \tfrac83\wp_{13}\wp_{1333}
  - 2{\wp_{12}}^2 +16\wp_{33}{\wp_{13}}^2  \nonumber \\ 
  & \quad + \mu_3(2\wp_{22}\wp_{12} - 2\wp_{33}\wp_{23}\wp_{13} 
  + \tfrac13\wp_{23}\wp_{1333})
  + (8\mu_6 - 2{\mu_3}^2)\wp_{33}\wp_{13}  \nonumber \\ & 
  \quad 
  - 2\mu_9 \wp_{33}\wp_{23},                         \quad [-14] \\
\wp_{333}\wp_{112} & = -2\wp_{23}{\wp_{13}}^2 + 2\mu_6\wp_{23}\wp_{13} 
  + 2\wp_{33}\wp_{23}\wp_{11} - 2\mu_{12}\wp_{23}
  - \mu_9{\wp_{23}}^2\nonumber \\ 
  & \quad+ \tfrac43 \wp_{12}\wp_{1333} - 4\mu_3{\wp_{13}}^2  
   - 4\wp_{33}\wp_{13}\wp_{12},                      \quad [-15]  \\
\wp_{113} \wp_{233} & = 2 \wp_{23}{\wp_{13}}^2
  + 2\wp_{33}\wp_{23}\wp_{11}
  - 2\mu_{12}\wp_{23} - \tfrac23\wp_{12}\wp_{1333} 
  + 2\mu_3{\wp_{13}}^2 \nonumber \\ 
  & \quad + 4\wp_{33}\wp_{13}\wp_{12} - \mu_9\wp_{33}\wp_{22} 
  + 2\mu_6\wp_{33}\wp_{12} + 2\mu_9\wp_{13},         \quad [-15]  \\
\wp_{123}\wp_{133} & = 2\wp_{23}{\wp_{13}}^2  
  - 2\mu_{12}\wp_{23} + \tfrac13 \wp_{12}\wp_{1333}
  + 2\mu_3{\wp_{13}}^2 + 2\mu_9\wp_{13},             \quad [-15]  \\
\wp_{122}\wp_{222} & = - \tfrac43\wp_{12}\wp_{1333} 
  + 8\wp_{33}\wp_{13}\wp_{12} + 4{\wp_{22}}^2\wp_{12}
  - 2\mu_3{\wp_{13}}^2 + 2\mu_3\wp_{33}\wp_{11}
  - \tfrac23 \mu_3\wp_{22}\wp_{1333}  \nonumber \\ 
  & \quad + 4\mu_3\wp_{33}\wp_{13}\wp_{22} 
  + {\mu_3}^2\wp_{23}\wp_{13}
  + (-2{\mu_3}^2+4\mu_6)\wp_{33}\wp_{12} 
  + (2 \mu_3 \mu_6-8\mu_9)\wp_{13} \nonumber \\ & 
  \quad + 6\mu_9\wp_{33}{\wp_{22}}
  - 2 \mu_9 {\wp_{23}}^2 -\mu_3 \mu_9 \wp_{23} 
  - 2\mu_3 \mu_{12}.                                 \quad [-15]
\end{align*}
The above equations describe the Jacobi variety as an algebraic
  variety, see also \cite{bel00} where a general matrix construction
  is given.  By eliminating odd powers with the aid of identities such
  as
\[
{\wp_{333}}^2\,{\wp_{233}}^2-(\wp_{333}\wp_{233})^2 = 0,
\]
one can obtain equations of the Kummer variety, $J/(u\rightarrow [-1]u)$.


\parindent=0pt
\begin{thebibliography}{10}
\addcontentsline{toc}{section}{References}
%Kopka~H, Daly~PW. 2003. \emph{A Guide to \LaTeX} (4th~edn).
%Addison-Wesley.
%
%Lamport~L. 1994. \emph{\LaTeX: a Document Preparation System} (2nd~edn).
%Addison-Wesley.%
%
%Mittelbach~F, Goossens~M. 2004. \emph{The \LaTeX\ Companion}
%(2nd~edn). Addison-Wesley.

\bibitem{ba97}
H.~F. Baker.
\newblock {\em Abelian Functions}.
\newblock Cambridge Univ. Press, Cambridge, 1897.

\bibitem{ba98}
H.~F. Baker.
\newblock On the hyperelliptic sigma functions.
\newblock {\em Amer. J. of Math.}, 20:301--384, 1898.


\bibitem{ba02}
H.~F. Baker.
\newblock On a system of differential equations leading to periodic function
\newblock {\em Acta Math.}, 26:135--156, 1902.


\bibitem{ba07}
H.~F. Baker.
\newblock {\em Multiply Periodic Functions}.
\newblock Cambridge Univ. Press, Cambridge, 1907.

\bibitem{bg06}
Sadie Baldwin and John Gibbons.
\newblock Genus 4 trigonal reduction of the {B}enney equations.
\newblock {\em J. Phys. A}, 39:3607--3639, 2006.

\bibitem{bego06}
Sadie Baldwin, J. C. Eilbeck, John Gibbons, and Y.~\^Onishi.
\newblock Abelian functions for purely trigonal curves of genus four.
\newblock {\tt http://arxiv.org/abs/math.AG/0612654},  2006.

\bibitem{bc28}
J.~L. Burchnall and T.~W. Chaundy.
\newblock Commutative ordinary differential operators. 
\newblock {\em Proc. London Math. Soc}, 118:420--440, 1923.

%\bibitem{bc28a}
%J.~L. Burchnall and T.~W. Chaundy.
%\newblock Commutative ordinary differential operators.
%\newblock {\em  Proc. Royal Soc. London (A)}, 118:557--583, 1928.

\bibitem{bel97}
V.~M. Buchstaber, V.~Z. Enolskii, and D.~V. Leykin.
\newblock Kleinian functions, hyperelliptic {J}acobians and applications.
\newblock {\em Reviews in Math. and Math. Physics}, 10:1--125, 1997.

%\bibitem{bel99}
%V.~M. Buchstaber, V.~Z. Enolskii, and D.~V. Leykin.
%\newblock Rational analogs of {A}belian functions.
%\newblock {\em Functional Anal. Appl.}, 33:83--94, 1999.

\bibitem{bel00}
V.~M. Buchstaber, V.~Z. Enolskii, and D.~V. Leykin.
\newblock Uniformization of {J}acobi varieties of trigonal curves 
and nonlinear  equations.
\newblock {\em Functional Anal. Appl.}, 34:159--171, 2000.

\bibitem{bl02}
V.~M. Buchstaber and D.~V. Leykin.
{Polynomial {L}ie algebras.} Func. Anal. Appl. 36:4, 267--280, 2002.

\bibitem{bl05}
V.~M. Buchstaber and D.~V. Leykin. {Addition Laws on Jacobian
Varieties of Plane Algebraic Curves} {\it Proceedings of the Steklov
Institute of Mathematics}, 251, 54-126, 2005

\bibitem{cn06}
K.~Cho and A.~Nakayashiki.
\newblock Differential structure of {A}belian functions.\\
\newblock {\tt http://arxiv.org/abs/math.AG/0604267}, 2006, to appear
in {\em International Journal of Mathematics}.

\bibitem{Weier34}
\newblock Weierstrass functions for higher genus curves. \hfill\newline
\newblock{\tt http://www.ma.hw.ac.uk/Weierstrass/}, maintained by
J.~C.~Eilbeck.

\bibitem{ee00} J.~C. Eilbeck, V.~Z. Enolskii.  
\newblock Bilinear operators and the power series for the Weierstrass $\sigma$
  function.
\newblock J. Phys. A. 33:791--794, 2000.

\bibitem{eel00} J.~C. Eilbeck, V.~Z. Enolskii, and D.~V. Leykin.
  \newblock On the {K}leinian construction of {A}belian functions of
  canonical algebraic curves.  \newblock In D.~Levi and O.~Ragnisco,
  editors, {\em Proceedings of the 1998 SIDE III Conference, 1998:
    Symmetries of Integrable Differences Equations}, volume CRMP/25 of
  {\em CRM Proceedings and Lecture Notes}, pages 121--138, 2000.

\bibitem{eep03}
J. C. Eilbeck, V. Z. Enolskii, and E. Previato.
\newblock On a generalized {F}robenius--{S}tickelberger addition formula.
\newblock Lett. Math. Phys, 65:5--17, 2003.

\bibitem{cafl96}
J.~W.~S.  Cassels and E. V. Flynn.
\newblock {\em Prolegomena to a middlebrow arithmetic of curves of genus $2$}
{\em London Math. Soc. Lect. Notes}, volume 230, Cambridge Univ. Press, 1996.

\bibitem{fa73}
J.~D. Fay, {\em Theta Functions on Riemann Surfaces}, Lecture Notes 
in Mathematics, vol.352, Springer-Verlag, Berlin, 1973.


\bibitem{gr90}
D.~Grant.
\newblock Formal groups in genus two.
\newblock {\em J. reine angew. Math.}, 411:96--121, 1990.

%\bibitem{ha84}
%G. H. Halphen.
%\newblock Memoire sur la reduction des equations differentielles
%lineaires aux formes integrales.
%\newblock {\em Mem. pres. l'Acade. de Sci. de France}, 28:1--300 (1884).

\bibitem{la82}
S.~Lang.
\newblock {\em Introduction to algebraic functions and {A}belian functions}.
\newblock Number~89 in Grad.Text in Math. Springer-Verlag, 2nd. edition, 1982.

\bibitem{mp05} S. Matsutani and E. Previato
\newblock {\em A Generalized Kiepert Formula for Plane Affine Curves}.
\newblock Institut Mittag-Leffler preprint, 2005 fall no. 03, 
http://www.mittag-leffler.se/preprints/ (ISSN 1103-467X ISRN 
IM-L-R-03-05/06--SE+fall).

\bibitem{mu83}
D.~Mumford.
\newblock {\em Tata Lectures on Theta I}, Progress in Mathematics, v. 28.
\newblock Birkh\"auser, 1983.

\bibitem{mu85}
D.~Mumford.
\newblock {\em Abelian varieties}.
\newblock Oxford Univ. Press, 1985.

\bibitem{on98}
Y.~\^Onishi.
\newblock Complex multiplication formulae for hyperelliptic curves of genus
  three.
\newblock {\em Tokyo J. Math.}, 21:381--431, 1998.
\newblock A list of corrections is available from {\tt
  http://web.cc.iwate-u.ac.jp/$\sim${}onishi/}

\bibitem{onPEMS} 
Y.~\^Onishi.  \newblock Determinant expressions for
  hyperelliptic curves (with an appendix by S.Matsutani) \newblock
  {\em Proc. Edinb. Math. Soc.}, 48:705--742, 2005.

\bibitem{on05}
Y.~\^Onishi.
\newblock Abelian functions for trigonal curves of degree four and
  determinantal formulae in purely trigonal case.
\newblock {\tt http://arxiv.org/abs/math.NT/0503696}, 2005.

\bibitem{shiga88}
H. Shiga.
\newblock On the representation of the {P}icard modular function
by $\theta$ constants {I}-{II} 
{\em Publ. RIMS, Kyoto Univ.}, 24:311-360, 1988.

\bibitem{sc89} 
R. J. Schilling.
\newblock Generalizations of the
Neumann system. A curve theoretical approach. II.
\newblock {\em Comm. Pure Appl. Math.},  42:409-442, 1989.

\bibitem{smb03}
Wolfgang Schreiner, Christian Mittermaier, and Karoly Bosa.
\newblock Distributed {M}aple: Parallel computer algebra in networked
  environments.
\newblock {\em Journal of Symbolic Computation}, 35:305--347, 2003.

\end{thebibliography}
\end{document}